\input amssym.def
\input amssym.tex

\headline={\ifnum \pageno=1 {\hfill} \else {\hss \tenrm -- \folio\ -- \hss}\fi}
\footline={\hfil}

\def\dater{\vglue-10mm\rightline{(\the\day/\the\month/\the\year)}}

\hsize 146mm
\vsize 224mm
\hoffset=6mm
\voffset=8mm
\baselineskip=5mm
\overfullrule =0pt

 at 10,5pt

\font\GGtitre=cmbx10 at 16pt


\def\og{\leavevmode\raise.30ex\hbox{$\scriptscriptstyle\langle\!\langle\>$}} 
\def\fg{\leavevmode\raise.24ex\hbox{$\scriptscriptstyle\>\rangle\!\rangle$}} 

\catcode`\@=11

\font\author=cmcsc10
\font\pauthor=cmcsc10 at 8pt
\font\tenmsx=msam10
\font\sevenmsx=msam10 scaled 700
\font\fivemsx=msam10 scaled 500
\font\tenmsy=msbm10
\font\sevenmsy=msbm10 scaled 700
\font\fivemsy=msbm10 scaled 500
\newfam\msxfam
\newfam\msyfam
\textfont\msxfam=\tenmsx  \scriptfont\msxfam=\sevenmsx
\scriptscriptfont\msxfam=\fivemsx
\textfont\msyfam=\tenmsy  \scriptfont\msyfam=\sevenmsy
\scriptscriptfont\msyfam=\fivemsy

\def\hexnumber@#1{\ifnum#1<10 \number#1\else
\ifnum#1=10 A\else\ifnum#1=11 B\else\ifnum#1=12 C\else
\ifnum#1=13 D\else\ifnum#1=14 E\else\ifnum#1=15 F\fi\fi\fi\fi\fi\fi\fi}

\def\msx@{\hexnumber@\msxfam}
\def\msy@{\hexnumber@\msyfam}
\mathchardef\nmid="3\msy@2D
\mathchardef\varnothing="0\msy@3F
\mathchardef\nexists="0\msy@40
\mathchardef\smallsetminus="2\msy@72
\def\Bbb{\ifmmode\let\next\Bbb@\else
\def\next{\errmessage{Use \string\Bbb\space only in math mode}}\fi\next}
\def\Bbb@#1{{\Bbb@@{#1}}}
\def\Bbb@@#1{\fam\msyfam#1}

\font\tentbl=cmr10 scaled 900
\font\seventbl=cmr7 scaled 900
\font\fivetbl=cmr5 scaled 900

\newfam\tblfam

\textfont\tblfam=\tentbl
\scriptfont\tblfam=\seventbl
\scriptscriptfont\tblfam=\fivetbl


\def \C {{\Bbb C}}

\def \F {{\Bbb F}}
\def \G {{\Bbb G}}
\def \h {{\Bbb H}}

\def \N {{\Bbb N}}
\def \P {{\Bbb P}}
\def \Q {{\Bbb Q}}
\def \R {{\Bbb R}}

\def \Z {{\Bbb Z}}

\font\ci= eufm10
\font\cisi= eufm9
\font\cis= eufm7
\font\ciss= eufm6
\font\cisss= eufm5

\def\ciB{\hbox{\ci B}}

\def\ciP{\hbox{\ci P}}

\def\cisiB{\hbox{\cisi B}}

\def\cisB{\hbox{\cis B}}

\def\cisP{\hbox{\cis P}}

\def\cissP{\hbox{\ciss P}}

\def\cisssP{\hbox{\cisss P}}

\font\ccm=cmmi10 at 9pt
\def\ccmL{\hbox{\ccm L}}

 at 9,5pt

\font\f=cmr10 at 7pt

\def\f1{\hbox{\f 1}}
\def\f2{\hbox{\f 2}}
\def\f3{\hbox{\f 3}}
\def\f4{\hbox{\f 4}}
\def\f5{\hbox{\f 5}}
\def\f6{\hbox{\f 6}}
\def\f7{\hbox{\f 7}}
\def\f8{\hbox{\f 8}}
\def\f9{\hbox{\f 9}}

\def \d {\,{\rm d}}

\def \dm {{\hbox {${1\over 2}$}}}
\def \sset {{\smallsetminus }}
\def\frac#1#2{ {{#1}\over{#2}}}
\def\eps{\varepsilon}
\def\Gal{\hbox{{\rm Gal}}}
\def\syms{\hbox{{\rm Sym}}^2\,}
\def\symm#1{\hbox{{\rm Sym}}^{#1}\,}

\def\Tr{\hbox{{\rm Tr}}\,}
\font\eutitre=eufm10 at 16pt
\font\cmtitre=cmmi10 at 16pt
\def\titreB{\hbox{\eutitre B}}
\def\titreL{\hbox{\cmtitre L}}
\def\sumb{\mathop{\sum \Bigl.^{\flat}}\limits}
\def\frob{\sigma}
\def\stacksum#1#2{{{\scriptstyle #1}}\atop {{\scriptstyle #2}}}

\topskip=10pt
\font\sept=cmti9

\def\rightheadline{\ifnum\pageno=\chstart{\hfill}
             \else{\centerline{\sept Gaps in coefficients of
$\ccmL$-functions and $\cisiB$-free numbers}}
             \hfill \hskip -30mm \tenrm\folio\fi}
\def\leftheadline{\ifnum\pageno=\chstart{\hfill}
             \else\tenrm\folio \hskip -3,5mm 
\hfill{\centerline{\pauthor E. Kowalski, O. Robert \& J. Wu}}\fi}
\headline={\ifnum\pageno=\chstart{\hfill}\else{\ifodd\pageno\rightheadline\else\leftheadline\fi}\fi}
\footline={\hfill}

\pageno=1
\newcount\chstart
\chstart=\pageno

\dater

\vglue 5mm

\centerline{\GGtitre Small gaps in coefficients of $\titreL$-functions}
\bigskip
\centerline{\GGtitre and $\titreB$-free numbers in short intervals}
\vskip 5mm

\centerline{\author E. Kowalski, O. Robert \& J. Wu}

\vskip 8mm

{\leftskip=1cm
\rightskip=1cm
{{\bf Abstract}.
We discuss questions related to the non-existence of gaps in
the series defining modular forms and other arithmetic functions of
various types, and improve results of Serre, Balog \& Ono and Alkan
using new results about exponential sums and the distribution of
$\ciB$-free numbers.
\par}}
\vskip 10mm

\noindent{\bf \S\ 1. Introduction}
\medskip

The motivation of this paper is a result of Serre
([43, Th. 15]) and the questions he subsequently raises. 
Let $f$
be a primitive holomorphic cusp form
(i.e. a newform in the Atkin-Lehner terminology) of weight $k$, with
conductor $N$ and nebentypus $\chi$. Write
$$
f(z)=\sum_{n\geq 1}{\lambda_f(n)e(nz)}
\leqno{(1.1)}
$$
its Fourier expansion at infinity, where $e(z)=\exp(2\pi iz)$, so that
$\lambda_f(n)$ is also the Hecke eigenvalue of $f$ for the Hecke
operator $T_n$. Serre's result is that  
$$
|\{p\leq x\,\mid\,\lambda_f(p)=0\}|\ll x(\log x)^{-1-\delta},
\leqno{(1.2)}
$$
for $x\geq 2$ and any $\delta<\frac{1}{2}$, the implied constant depending
on $f$ and $\delta$, from which he deduces that the
series~(1.1), or equivalently the $L$-function
$$
L(f,s)=\sum_{n\geq 1}{\lambda_f(n)n^{-s}}
\leqno{(1.3)}
$$
is not lacunary, i.e. the set of indices $n$ where $\lambda_f(n)\not=0$,
has a positive density. Serre asked ([43, p. 183]) for
more precise statements, in particular for bounding non-trivially the
function $i_f(n)$ defined by
$$
i_f(n)=\max\{ k\geq 1\,\mid\,\lambda_f(n+j)=0\hbox{ for } 0<j\leq k\},
\leqno{(1.4)}
$$
where non-trivial means an estimate of type $i_f(n)\ll n^{\theta}$ for
some $\theta<1$ and all $n\geq 1$.
A stronger form of the problem is to find $y$ as small as possible (as
a function of $x$, say $y=x^{\theta}$ with $\theta<1$) such that
$$
|\{n\,\mid\, x<n\leq x+y\hbox{ and } \lambda_f(n)\not=0\}|\gg y
\leqno{(1.5)}
$$
(where the implied constant can depend on $f$). Non-lacunarity means
$y=x$ is permitted, and one wishes to improve this. Note $i_f(n)\ll y$
so this generalizes the first question. 
\par
\medskip
The history of this problem is somewhat confused. First, Serre could
have solved it quite simply in (at least) two ways available at the
time. The first is to argue that by multiplicativity
$\lambda_f(n)\not=0$ if $n$ is squarefree and not divisible by primes
$p$ for which $\lambda_f(p)=0$. The latter have density zero
by~(1.1), so estimating $i_f(n)$ becomes a special case of
a problem in multiplicative number theory, that of counting so-called 
$\ciB$-free numbers in small intervals, where for a set
$\ciB=\{b_i\}$ of integers with $(b_i,b_j)=1$ if $i\not=j$ and
$$
\sum_{i}{\frac{1}{b_i}}<+\infty,
$$
one says that $n\geq 1$ is $\ciB$-free if it is not divisible
by any element in $\ciB$. Erd\"os [11] already showed in 1966
that with no further condition there exists a constant $\theta<1$
(absolute) such that the interval $(x,x+x^{\theta}]$
contains a 
$\ciB$-free number for $x$ large enough, thereby solving
Serre's first question in the affirmative. A quantitative result
proving the analogue of~(1.5) for general $\ciB$-free
numbers was also obtained Szemer\'edi [44] as early as 1973.
\par
This was apparently first noticed by Balog and Ono [2]. By this time
results about $\ciB$-free numbers had 
been refined a number of times, and they deduced from a result of Wu
[45] that $i_f(n)\ll n^{17/41+\eps}$ for $n\geq 1$ and any $\eps>0$, the
implied constant depending on $f$ and $\eps$. Using this idea and
other results (such as a version of the Chebotarev density theorem in
small intervals and the Shimura correspondence), they also get
weaker results for modular forms of weight $1$
or half-integral weight. The latter is noteworthy in this respect
since the Fourier coefficients of half-integral weight forms are
highly non-multiplicative (see [7] for a strong
quantitative expression of this fact).
Alkan [1] has developed and
improved the results of [2], tailoring some arguments to the specific
instance of $\ciB$-free numbers involved for the problem at
hand. 
\par
A second method of estimating $i_f(n)$ available to Serre was a direct
appeal to the properties of the Rankin-Selberg $L$-function
$L(f\otimes \bar{f},s)$.  Specifically this proves [36, 42] (for $f$ any
cusp form of integral weight $k$ and level $N$) 
$$
\sum_{n\leq x}{|\lambda_f(n)|^2n^{1-k}}=c_f x+O(x^{3/5})
$$
for some $c_f>0$, and $x\geq 1$, the implied constant
depending only on $f$. Trivially this implies $i_f(n)\ll
n^{3/5}$, and incidentally this fact is implicit in [27] (which
Serre quotes as one source for his problems!)
\par
It turns out however that there are still a number of things which
seem to have been overlooked. For instance we will show that it is
not necessary to sieve by squarefree numbers, and we will explain the
applications of the Rankin-Selberg $L$-functions (in particular to
non-congruence subgroups, another of the questions in [43]).
We also look at lacunarity in some other Dirichlet series coming from
arithmetic or analysis, including one
which is really neither fish nor fowl (see Proposition 4).
On the other hand (this is our main new contribution), we
will improve quite significantly the 
$\ciB$-free number results that can be used. Some of our tools are new
estimates for exponential sums and bilinear forms which are of
independent interest in analytic number theory.
\par
We of course welcome any further corrections to the picture thus produced
about this problem.
\par
\bigskip

\noindent{\bf Acknowledgement.}
The authors would like to thank Emmanuel Royer
for helpful comments on an earlier version of this paper.
\par
\medskip
\noindent{\bf Notation.} For any $k\geq 1$, $N\geq 1$ and any
character $\chi$ modulo $N$, we denote $S_k(N,\chi)$ the vector space
of cusp forms of weight $k$ for the group $\Gamma_0(N)$, with
nebentypus $\chi$. If $\chi$ is the trivial character modulo $N$, we
simply write $S_k(N)$. We also denote by $S_k^*(N,\chi)$, or
$S_k^*(N)$, the set of primitive forms in $S_k(N,\chi)$ or $S_k(N)$,
i.e. those forms which are eigenfunctions of all Hecke operators $T_n$
and are normalized by $\lambda_f(1)=1$, where $\lambda_f(n)$ is the
$n$-th Fourier coefficient, which is then equal to the $n$-th Hecke
eigenvalue. See e.g. [23] for basic analytic facts about modular
forms. 
\par
For $s$ a complex number, we denote $\sigma$ its real part and
$t$ its imaginary part. Also, we use $f(x)=O(g(x))$ and $f(x)\ll g(x)$
for $x$ in some set $X$ as synonyms, meaning $|f(x)|\leq Cg(x)$ for
all $x\in X$, $C\geq 0$ being called the implied constant.

\vskip 5mm
\goodbreak
\noindent{\bf \S\ 2. Algebraic aspects}
\medskip

We start by noticing that the restriction to squarefree numbers
present in [2] and [1] is in fact unnecessary, because the set of
primes for which $\lambda_f(p^{\nu})=0$ for {\it any} $\nu$ still
satisfies an estimate similar to (1.2). This is partly implicit in
[43, p. 178--179]. 

\proclaim Lemma 2.1. Let $f\in S_k^*(N,\chi)$ be a primitive holomorphic
cusp form. There exists an integer $\nu_f$ depending only on
$f$ such that for any prime $p\,\nmid N$, either
$\lambda_f(p^{\nu})\not=0$ for all $\nu\geq 0$, or there exists
$\nu\leq \nu_f$ such that $\lambda_f(p^{\nu})=0$.

\noindent{\sl Proof}.
Let $p\,\nmid N$.
By multiplicativity we have the power series expansion
$$\sum_{\nu\ge 0} \lambda_f(p^\nu) X^\nu
= {1\over 1-\lambda_f(p)X+\chi(p)p^{k-1}X^2}.
\leqno(2.1)$$
Let $\alpha_p$ and $\beta_p$ be the complex numbers such that
$$1-\lambda_f(p)X+\chi(p)p^{k-1}X^2 = (1-\alpha_p X) (1-\beta_p X).
\leqno(2.2)$$
Thus
$$\alpha_p+\beta_p = \lambda_f(p)
\qquad{\rm and}\qquad
\alpha_p\beta_p = \chi(p)p^{k-1}\not=0.
\leqno(2.3)$$
Expanding (2.1) using (2.2) by geometric series gives the
well-known expressions
$$\lambda_f(p^\nu)
= {\alpha_p^{\nu+1}-\beta_p^{\nu+1}\over \alpha_p-\beta_p}
\leqno(2.4)$$
if $\alpha_p\not=\beta_p$, and the simpler
$$\lambda_f(p^\nu)=(\nu+1)\alpha_p^\nu=(\nu+1) \tau
p^{\nu(k-1)/2}\not=0
\leqno(2.5)$$
if $\alpha_p=\beta_p$, where $\tau^2=\chi(p)$. So we can assume
$\alpha_p\not=\beta_p$. 
In this case we get by (2.4)
$$\lambda_f(p^\nu)=0
\quad\hbox{if and only if}\quad
(\alpha_p/\beta_p)^{\nu+1}=1
$$
so that there exists $\nu\ge 0$ for which $\lambda_f(p^\nu)=0$ if and
only if $\alpha_p/\beta_p$ is a root of unity, and if this ratio is a
primitive root of unity of order $d\ge 1$, then $\lambda_f(p^{d-1})=0$.
\par
Now we input some more algebraic properties of the Fourier
coefficients. The field
$$
K_f = \Q(\lambda_f(n),\chi(n))
$$
generated by all Fourier coefficients and values of $\chi$ is known to
be a number field. By (2.2), the ``roots'' $\alpha_p$ and
$\beta_p$ lie in a quadratic extension of $K_f$. This extension (say
$K_p$) depends on $p$, but it has degree $[K_p : \Q]\le 2[K_f : \Q]$
for all $p$.
\par
Now we combine both remarks and the fact that a number field $L/\Q$ can
only contain a primitive $d$-th root of unity if $\varphi(d)\le [ L: \Q]$.
It follows that if $p\,\nmid N$ and $\lambda_f(p^\nu)=0$ for some
$\nu\ge 0$, $\alpha_p/\beta_p$ is a primitive root of unity of some
order $d$ such that $\varphi(d)\le 2[K_f : \Q]$, and then
$\lambda_f(p^{d-1})=0$. Since $\varphi(d)\gg d/\log\log(3d)$,
this proves the lemma.\hfill$\square$
\smallskip

It is clear that $\nu_f$ is effectively computable. Here are some
simple cases.

\proclaim Lemma 2.2. Let $k$ be even and $f\in  S_k^*(N,\chi)$. There
exists $M\geq 1$ such that for any $p\,\nmid M$, either
$\lambda_f(p)=0$ or $\lambda_f(p^{\nu})\not=0$ for any $\nu\geq 1$. If
$\chi$ is trivial and $f$ has integer coefficients, one can take $M=N$.

\noindent{\sl Proof}.
If $p\mid N$, the condition $\lambda_f(p^{\nu})=0$ is equivalent to
$\lambda_f(p)=0$ by total multiplicativity, so we can assume that
$p\,\nmid N$. Let $p$ be such a prime with $\lambda_f(p^\nu)=0$ for
some $\nu\ge 2$, 
but $\lambda_f(p)\not=0$. Using the same notation as the proof of
Lemma 2.1, we have $\alpha_p=\xi \beta_p$ for some root of
unity $\xi$ of order $d+1$, and $\xi\not=-1$.
We derive from the second relation of (2.3) that
$\alpha_p^2 = \xi \chi(p)p^{k-1}$,
hence $\alpha_p = \pm \tau p^{(k-1)/2}$,
where $\tau^2 = \xi \chi(p)$.
By the second relation of (2.3), we get
$$\lambda_f(p) = (1+\bar{\xi})\alpha_p=\pm
\tau(1+\bar{\xi})p^{(k-1)/2}\not=0.$$ 
In particular, since $k$ is even,
$\Q(\tau(1+\bar{\xi})\sqrt{p})\subset K_f$.
As $K_f$ is a number field,
this can happen only for finitely many $p$, and one can take as $M$
the product of those primes and those $p\mid N$ with $\lambda_f(p)=0$.
\par
Furthermore, if $\chi$ is trivial and $f$ has integer coefficients,
then for $p\nmid N$, $\alpha_p/\beta_p=\xi$ is a root of unity
$\not=1$ in a quadratic  
extension of $\Q$ (see (2.2)),
hence $\xi\in\{ -1, \pm j,\pm j^2,\pm i\}$ (with $\nu\in \{1, 2, 3, 5\}$).
All those except $\xi=-1$ contradict the fact that $f$ has 
integer coefficients by simple considerations such
as the following, for $\xi=j$ say:
we have $\alpha_p^2 = j p^{k-1}$, $\alpha_p = \pm j^2 
p^{(k-1)/2}$ and
$$\lambda_f(p)
= (1+\bar{j})\alpha_p
= \pm (1+\bar{j})j^2p^{(k-1)/2}
= \pm (j^2+j)p^{(k-1)/2}\notin\Z$$
(compare [43, p. 178--179]).
\hfill$\square$
\smallskip

We now prove the analogue of (1.2) for primes $p$ such
that $\lambda_f(p^{\nu})=0$ for some $\nu$.

\proclaim Lemma 2.3. Let $f\in S_k^*(N,\chi)$ be a primitive cusp form
not of CM type, in particular with $k\geq 2$. For $\nu\geq 1$, let
$$
\ciP_{f, \nu}=\{p\,\nmid N\,\mid\, \lambda_f(p^{\nu})=0\}.
\leqno{(2.6)}
$$
For any $\nu\geq 1$ we have
$$\big|\ciP_{f, \nu} \cap [1, x]\big|
\ll {x\over (\log x)^{1+\delta}}
\leqno(2.7)$$
for $x\ge 2$ and any $\delta<\dm$, the implied constant depending on
$f$ and $\delta$. Let $\ciP_{f}^*$ be the union of $\ciP_{f,\nu}$. We
have 
$$
\big|\ciP_{f}^* \cap [1, x]\big|
\ll {x\over (\log x)^{1+\delta}}
\leqno(2.8)
$$
for $x\geq 2$ and any $\delta<\dm$, the implied constant depending on
$f$ and $\delta$.

\noindent{\sl Proof}.
All the tools needed to prove (2.7), if not the exact statements,
can be gathered from [43], in particular Section 7.2.
By Lemma 2.1, we need only prove (2.7), so let $\nu\geq 1$ be fixed.
\par
Fix a prime number $\ell$ totally split in the field
$K_f=\Q(\lambda_f(n),\chi(n))$ already considered. Thus $K_f\subset
\Q_{\ell}$. There exists an
$\ell$-adic Galois representation
$$
\rho_{f,\ell} \,:\, {\rm Gal}(\bar{\Q}/\Q)\to GL(2,\Q_{\ell})
$$
constructed by Deligne, such that for $p\,\nmid N\ell$ we have
$${\rm Tr}\rho_{f,\ell}(\sigma_p) = \lambda_f(p)
\qquad{\rm and}\qquad
\det\rho_{f,\ell}(\sigma_p)=\chi(p)p^{k-1},$$
where $\sigma_p$ is a Frobenius at $p$. Let $G_{\ell}$ be the image 
of $\rho_{f,\ell}$.
As explained by Serre [43, Prop. 17], it is an open
subgroup of $GL(2,\Q_{\ell})$, hence an $\ell$-adic group of
dimension $4$.
\par
By symmetry, there exists a polynomial $P_\nu\in \Z[X,Y]$ such that
the identity
$${X^{\nu+1}-Y^{\nu+1}\over X-Y} = P_\nu(X+Y, XY)$$
holds.
Consider the set $C\subset G_{\ell}$ defined by
$$C = \{s\in G_{\ell}\,\mid\, P_\nu({\rm Tr}(s), \det(s))=0\}.$$
Note the following facts about $C$: it is a closed $\ell$-adic
subvariety of $G_{\ell}$, stable by conjugation, and of
dimension $\le 3$. Moreover, $C$ is stable by multiplication by
$H_{\ell}=\{\hbox{homotheties in $G_{\ell}$}\}$, and therefore
$C=\pi^{-1}(C')$ for a certain subvariety $C'\subset
G_{\ell}/H_{\ell}$,
where $\pi : \, G_{\ell}\to G_{\ell}/H_{\ell}$ is the projection.
The set $C'$ is an $\ell$-adic variety of
dimension $\le 2$ and all its elements are regular ([43, Section 5.2]),
since they have distinct eigenvalues $\alpha$,
$\xi\alpha$ for some root of unity $\xi\not=1$ of order $\nu+1$.
\par
Now remark that if $p\in \ciP_{f, \nu}$ and $p\,\nmid N\ell$, we have
$\pi(\sigma_p)\in C'$ (going back to the proof Lemma 2.1 if necessary).
Hence our result (2.7) follows from Theorem 12 of [43],
as in the proof of the case $h=0$ of Theorem 15 of
loc. cit., p. 177.
\hfill
$\square$
\smallskip

For ease of reference we recall the lemma which allows the extension
of the results for $i_f(n)$ to general cusp forms from that of
newforms.

\proclaim Lemma 2.4. Let $f\in S_k(N,\chi)$ be a cusp form not in the
space spanned by CM forms. There exist:

\vskip -2mm
{\sl
{\rm (i)} an integer $s\geq 1$ and
algebraic numbers $\beta_j$ and
positive rational numbers $\gamma_j$ for $1\leq j\leq s$;
\par
{\rm (ii)} a divisor $\delta\mid N$ such that $\chi$ is induced by 
$\chi_1$ modulo $N/\delta$ and a divisor $\delta_1\mid \delta$;
\par
{\rm (iii)} a primitive form $g\in S_k^*(N/\delta,\chi_1)$, not
of CM type;
\par
\noindent
such that
$$
\lambda_g(n)=\sum_{1\leq j\leq s}{\beta_j\lambda_f(\gamma_j\delta_1 n)}
$$
for $n\geq 1$. By convention, we put $\lambda_f(x)=0$ if $x\in \Q$ is
not a positive integer.}

\smallskip
This is just a formal restatement of the computations in
[2], p. 362, or follows from [43, \S 7.6].
\par
\smallskip
We now discuss briefly the possibility of extending the results above
to higher rank situations. From the proof of Lemma 2.3, it is natural
to start from an $\ell$-adic representation
$$
\rho\,:\, \hbox{Gal}(\bar{\Q}/\Q)\rightarrow GL(V)
$$
where $V\simeq \Q_{\ell}^r$ for some $r\geq 1$. We assume it is
``sufficiently geometric'', namely that it is unramified outside a finite
set of primes $S$, and that the 
$L$-function of $\rho$, defined as usual by the Euler product
$$
L(\rho,s)=\prod_p{\det(1-\rho(\frob_p)p^{-s}\mid V^{I_p})^{-1}}
=\sum_{n\geq 1}{\lambda_{\rho}(n)n^{-s}},
\leqno{(2.9)}
$$
(where $\frob_p$ is a Frobenius element at $p$ and $I_p$ the inertia
group at $p$) has coefficients in a number field $K_{\rho}\subset
\Q_{\ell}$. Note that we view this here as a formal Dirichlet series.
If the image of $\rho$ is fairly big, one can use the methods of Serre
to get
$$
|\{p\leq x\,\mid\, p\notin S \hbox{ and }
\lambda_{\rho}(p)=\Tr\rho(\frob_p)=0\}|\ll x(\log x)^{-1-\delta}
$$
for some $\delta>0$, see Proposition 1 below.
On the other hand, it is not clear if the 
analogue of Lemma 2.1 holds, and this seems a hard question in
general. The analogue of~(2.4) does not provide an equation
easily solvable to characterize the values of $\nu$ for which
$\lambda_{\rho}(p^{\nu})=0$. The best that seems doable is to notice that,
for fixed (unramified) $p$, 
$u_{\nu}=\lambda_{\rho}(p^{\nu})$ is given by a linear recurrence relation of
degree $r$ with ``companion polynomial'' given by
$$
\det(X-\rho(\frob_p))=
X^r-\lambda_{\rho}(p)X^{r-1}+\cdots+ (-1)^r\det\rho(\frob_p)=
\prod_{1\leq i\leq r}{(X-\alpha_{p,i})}
$$
so that
$$
u_{\nu}=\lambda_{\rho}(p^{\nu})=\sum_{1\leq i\leq
r}{\gamma_{p,i}\alpha_{p,i}^{\nu}} 
$$
for some $\gamma_{p,i}$. 
The Skolem-Mahler-Lech theorem (see e.g. [6, p. 88]) says that
for any linear recurrence 
sequence $(u_{\nu})$, either $u_{\nu}=0$ for only finitely many values
of $\nu$, or there exists an arithmetic progression $a+tv$ with $v\not=0$ such
that $u_{a+tv}=0$ for all $t\geq 0$. In the latter case, spelling this
out yields a Vandermonde type linear system for powers of the
$\alpha_{p,i}^v$, hence it implies that
$\alpha_{p,i}^v=\alpha_{p,j}^v$ for some $i\not=j$. Coming back to
$\rho$, this case implies that an extension of degree $\leq
r[K_{\rho}:\Q]$ contains a $v$-th root of unity. As in
Lemma 2.1, this bounds $v$, and an analogue of
Lemma 2.3 is possible, given $\xi$ a root of unity, to get
$$
|\{p\leq x\,\mid\, p\notin S\hbox{ and there are two roots
  $\alpha_{i,p}$, $\alpha_{j,p}=\xi\alpha_{i,p}$ with $i\not=j$}
\}|\ll x(\log x)^{-1-\delta}
$$
for some $\delta>0$.
\par
However, in the first case where $u_{\nu}=0$ has only finitely many
solutions, despite the remarkable fact that there exists a uniform
bound for the number of solutions depending only on $r$ (see
[12]), this is insufficient because only the number of
solutions, not the value of $\nu$, is bounded, so that an integer $\nu_0$
(independent of $p$) for which the smallest solution is $\nu\leq \nu_0$ is
not known to exist. The question amounts to asking for a bound for the
height of the solutions to the relevant linear equations in multiplicative
groups [12, p. 820], and is thus in full generality of the same
type as asking for effective versions of Roth's theorem, or of
Schmidt's Subspace Theorem. (Note that by replacing $u_{\nu}$ by
$p^{-(k-1)/2}u_{\nu}$ one gets a linear recurrence relation with companion
polynomial having height absolutely bounded, by the
Ramanujan-Petersson conjecture proved by Deligne).
\par
The theory of $\ciB$-free numbers does however still
apply. Thus we get:

\proclaim Proposition 1.
Let $\rho$ be an $\ell$-adic representation of
$\Gal(\bar{\Q}/\Q)$ on $GL(r,\Q_{\ell})$. Assume that $\rho$ is
unramified for $p$ outside a finite set $S$ and that its $L$-function
has coefficients in a number field $K_{\rho}$. Let
$G=\hbox{{\rm Im}}\, \rho$, $C=G\cap \{s\in GL(r,\Q_{\ell})\,\mid\, \Tr
s=0\}$. Assume that, 
as $\ell$-adic varieties, we have $\dim C<\dim G$. Then 
for any $\eps>0$, $x\geq x_0(\rho,\eps)$ and $y\geq x^{7/17+\eps}$ we
have 
$$
|\{n\,\mid\, x<n\leq x+y\hbox{ and } \lambda_{\rho}(n)\not=0\}|\gg y.
$$
In particular $i_{\rho}(n)\ll n^{7/17+\eps}$. 

\noindent{\sl Proof}.
One can argue as for modular forms using $\ciB$-free numbers (see
Proposition 6) with
$$
\ciB=\{ p\,\mid\, p\in S \hbox{ or }\lambda_{\rho}(p)=0\}
\cup \{p^2\,\mid\,
\lambda_{\rho}(p)\not=0\}, 
$$
after applying Theorem 10 of [43] to $G$ and $C$, with
$E=\bar{\Q}^{\ker \rho}$, to derive
$$
|\{p\le x\,\mid\, p\notin S \hbox{ or } \lambda_\rho(p)=0\}|
\ll x(\log x)^{-1-\delta}
$$
for some $\delta>0$ depending on the dimensions
of $G$ and $C$ (for instance, any $\delta<1-\dim C/\dim G$). Strictly
speaking, to apply this theorem as stated we must also treat
separately the  case where $G$ is finite. One can then see $\rho$ as a
linear representation of the finite group $G=\Gal(E/\Q)$ into
$GL(n,\bar{\Q})$, or into $GL(n,\C)$. In that case the condition $\dim 
C<\dim G$ means that the character of $\rho$ does not vanish. By
a well-known fact about linear representations of finite groups
(see e.g. [17, Ex. 2.39]),
this means that the representation $\rho$ is a one-dimensional
character of $\Gal(\bar{\Q}/\Q)$, which by the Kronecker-Weber
theorem corresponds to (i.e. has the same $L$-function as) a Dirichlet
character $\chi$, of conductor $N$ say (divisible only by primes in
$S$). Then $\lambda_{\rho}(n)\not=0$ if and only if $(n,N)=1$. 
\hfill$\square$
\smallskip

This is also implicit in [43, \S 6.4, 6.5].

\vskip 5mm
\noindent{\bf \S\ 3. Maass forms, cofinite groups and the
Rankin-Selberg method} 
\medskip

In this section, we describe what results follow from the
Rankin-Selberg method. Although, for fixed $f\in S_k^*(N,\chi)$, they
are weaker than those obtained by means of $\ciB$-free numbers, this
method has the advantage of yielding quite easily estimates uniform in
terms of $f$, i.e. with explicit dependency on $k$ and $N$. Those are
by no means obvious from the $\ell$-adic point of view leading to
(1.2). Moreover, the Rankin-Selberg method applies, at least as far
as bounding $i_f(n)$, to non-congruence subgroups, as shown by Good
[19], Sarnak [41] and Petridis [34]. This answers the last question
in [43, p. 183].

\proclaim Proposition 2. Let $\Gamma\subset SL(2,\R)$ be a discrete
subgroup such that the quotient $\Gamma\backslash\h$ has finite
hyperbolic volume and $\Gamma$ contains the integral translation
matrices acting by $z\mapsto z+n$. Let $f$ be either a holomorphic
cusp form of weight $k\geq 2$ or a Maass cusp form with eigenvalue
$\lambda\not=1/4$. Define $i_f(n)$ by (1.4) where $\lambda_f(n)$ are
the Fourier coefficients in the expansion of $f$ at the cusp
$\infty$ of $\Gamma$. Then for some $\theta<1$ we have
$$
i_f(n)\ll n^{\theta}
$$
for $n\geq 1$ where the implied constant depends on
$f$. Specifically, one can take $\theta=2/3$ if $f$ is
holomorphic and any $\theta>4/5$ if $f$ is non-holomorphic.

\noindent{\sl Proof}. The non-holomorphic case follows from [34] as the
holomorphic case follows from [19], so we describe only the
latter. Good shows that
$$
\sum_{n\leq x}{|\lambda_f(n)|^2}=
\sum_{2/3<s_j\leq 1}{\frac{(4\pi x)^{s_j+k-1}}{\Gamma(k+s_j)}
\langle r_j,y^k|f|^2\rangle}
+O(x^{k-1+2/3})
\leqno{(3.1)}
$$
for $x\geq 1$, where $1=s_0>s_1\geq \cdots\geq s_r$ are the finitely
many poles of the Eisenstein series $E(z,s)$ for $\Gamma$ in the
interval $[1/2,1]$ (those with $s_j>2/3$ go to the error
term), $r_j(z)$ is the residue of $E(z,s)$ at $s_j$ and
$\langle \cdot,\cdot\rangle$ is the inner product on
$L^2(\Gamma\backslash\h)$. The pole at $s_0=1$ with residue $V^{-1}$
contributes 
$$
\frac{(4\pi x)^{k}\|f\|^2}{k!V}
$$
where $\|f\|$ is the Petersson norm of $f$ and $V$ the volume of
$\Gamma\backslash\h$. Comparing (3.1) at
$x=n$ and $x=n+Cn^{2/3}$, where 
$C$ is some large constant, shows that $i_f(n)\leq Cn^{2/3}$.
\hfill $\square$

\smallskip
{\sl Remark 1.} As for half-integral weight forms, it is not expected
that the coefficients of a cusp form for a non-arithmetic group
satisfy any multiplicativity properties. In fact, it would be quite
interesting to express this in a quantitative manner as done by Duke
and Iwaniec [7] for half-integral forms using bilinear forms in the
Fourier coefficients.
\par
\smallskip
In the case of congruence subgroup the methods using $\ciB$-free
numbers yield stronger results such as (1.5) for $y$ quite
small. Those however are not uniform in terms of $f$ (i.e. in terms of
$N$ and $k$ for holomorphic forms). The Rankin-Selberg method can
quite easily yield some uniform estimates. Here are sample statements;
note that we have not tried to get the best possible results.

\proclaim Proposition 3. {\rm (1)} Let $N\geq 1$ and $f$ a primitive Maass
form of conductor $N$ with eigenvalue $\lambda\not=1/4$ and trivial
nebentypus. Let $\Lambda=\lambda+1$ and let
$\lambda_f(n)$ be the Fourier coefficients of $f$. For any $\eps>0$,
there exists $c>0$ depending only on $\eps$ such that when
$$
y>x^{37/40}\Lambda^{45/32}N^{19/8}(x\Lambda N)^{\eps},
\leqno{(3.2)}
$$
we have
$$
|\{n\,\mid\, x<n\leq x+y\hbox{ and } |\lambda_f(n)|^2\geq c(\Lambda
N)^{-\eps}\}|\gg y (\log x)^{-1}(\Lambda N)^{-\eps},
$$
the implied constants depending only on $\eps$.
\par
\vskip -2mm
{\rm (2)} {\sl Let $N\geq 1$ and $f$ a primitive holomorphic
form of conductor $N$, weight $k$ with nebentypus $\chi$, not of CM
type. There exists an absolute constant $c>0$ such that for any
$\eps>0$ and 
$$
y>x^{4/5}(kN)^{1/2}(xkN)^{\eps}
$$
we have
$$
|\{n\,\mid\, x<n\leq x+y\hbox{ and } |\lambda_f(n)|^2n^{(1-k)}\geq c(\log 
kN)^{-1}\}|\gg y (\log x)^{-14}(\log kN)^{-3}
$$
the implied constant depending only on $\eps$.}

\noindent{\sl Proof}. We prove (1) and only give some indications for
the easier (2) at the end. It turns out to be simpler to reduce to
squarefree numbers (so in fact we could impose this condition on
$n$). The result will follow by Cauchy's inequality from the two
asymptotic formulas
$$
\sumb_{n\leq x}{|\lambda_f(n)|^2}=c_f x + 
O((\Lambda N)^{25/64}x^{121/146}(x\Lambda N)^{\eps})
\leqno{(3.3)}
$$
and
$$
\sumb_{n\leq x}{|\lambda_f(n)|^4}=d_f x\log x+ e_f x+
O(\Lambda^{45/32}N^{19/8}x^{37/40}(x\Lambda N)^{\eps})
\leqno{(3.4)}
$$
where $\sumb$ restricts $n$ to squarefree integers coprime with
$N$. Both hold for any $\eps>0$, with the implied constant
depending only on $\eps$, and $c_f$, $d_f$, $e_f$ are real numbers
with $c_f$, $d_f>0$ and 
$$
c_f\gg (\Lambda N)^{-\eps},\quad d_f,\ e_f\ll (\Lambda N)^{\eps}
\leqno{(3.5)}
$$
for any $\eps>0$, the implied constant depending only on $\eps$. 
\par
Indeed, let $\eps>0$ and $\eta>0$ be any positive numbers, and put the
integers $n$ 
with $x<n\leq x+y$ in two sets $L$ and $S$ if, respectively,
$|\lambda_f(n)|^2>\eta$ or $|\lambda_f(n)|^2\leq \eta$. If $y$ satisfies
(3.2) we have by (3.3) and (3.5)
$$
\sumb_{x<n\leq x+y}{|\lambda_f(n)|^2}\gg c_f y\geq C y (\Lambda N)^{-\eps/4}
$$
where $C$ depends only on $\eps$, whereas by positivity and Cauchy's
inequality 
$$
\sumb_{x<n\leq x+y}{|\lambda_f(n)|^2}
\leq \eta |S| + \sumb_{x\in L}{|\lambda_f(n)|^2}
\leq \eta y + |L|^{1/2} \Bigl(
\sumb_{x<n\leq x+y}{|\lambda_f(n)|^4}\Bigr)^{1/2}.
$$
If $\eta<\frac{C}{2}(\Lambda N)^{-\eps/4}$ we derive by (3.4) and (3.5)
$$
|L|\gg y^2(\Lambda N)^{-\eps}y^{-1}(\log x)^{-1},
$$
as desired.
\par
We give the proof of (3.4) and the upper bounds on $d_f$, $e_f$,
since (3.3) is easier. The lower bound for $c_f$ is deeper,
and follows immediately from the bound $L(F,1)\gg (\Lambda N)^{-\eps}$
of Hoffstein and Lockhart [22] for the adjoint square $F$ of $f$
(which is also its symmetric square since the nebentypus is trivial).
\par
Since $f$ is primitive and has trivial nebentypus, hence real
coefficients, we have 
$$
|\lambda_f(p)|^4=\lambda_f(p)^4=
(1+\lambda_f(p^2))^2=1+2\lambda_f(p^2)+
\lambda_f(p^2)^2
$$
for $p\,\nmid N$ and thus we find that
$$
\eqalign{
L(s)&:=\sumb_{n\geq 1}{|\lambda_f(n)|^4n^{-s}}=
\prod_{p\,\nmid N}{(1+|\lambda_f(p)|^4p^{-s})}\cr
&=\zeta^{\flat}(s)L^{\flat}(\syms f,s)^2L^{\flat}(\syms f\otimes
\syms f,s)H(s)}
$$
where 
$$
\eqalign{
\zeta^{\flat}(s)&=\prod_{p\,\nmid N}{(1+p^{-s})},
\cr
L^{\flat}(\syms f,s)&=\prod_{p\,\nmid N}{(1+\lambda_f(p^2)p^{-s})},\cr
L^{\flat}(\syms f\otimes \syms f,s)&=
\prod_{p\,\nmid N}{(1+\lambda_f(p^2)^2p^{-s})}
}
$$
and $H(s)$ is an Euler product which converges absolutely for
$\sigma>\frac{23}{32}$ by the estimate $|\lambda_f(p)|\leq 2p^{7/64}$
of Kim and Sarnak [26] (any estimate $|\lambda_f(p)|\leq 2p^{\theta}$
with $\theta<1/4$ would do, at the cost of worsening the exponent),
and is moreover uniformly bounded (in terms of $f$) on any line
$\sigma=\sigma_0$ with $\sigma_0>\frac{23}{32}$. To see this, define
$H$ as the obvious ratio for $\sigma$ large enough, and check on the
Euler factors individually that
$$
H(s)\ll \zeta(2\sigma-14/32)^B
$$
for some absolute constant $B>0$ (for a similar argument, see e.g. [9,
Prop. 2]). 
\par
Each of the three $L$-functions is obtained by removing non-squarefree
coefficients (and those not coprime with $N$) from an $L$-function
which has analytic continuation and a functional equation of the
standard type: the first one is the zeta function, the second one is
the adjoint square $F=\syms f$ of Shimura and Gelbart-Jacquet [18],
and the third is the Rankin-Selberg square $F\otimes F$ of the latter
(which exists as a special case of convolution of cusp forms on
$GL(3)$). The same bound and reasoning already used shows that
$$
\zeta^{\flat}(s)L^{\flat}(F,s)^2L^{\flat}(F\otimes F,s)=
\zeta(s)L(F,s)^2L(F\otimes F,s)H_1(s)
$$
where $H_1(s)$ has the same properties as $H(s)$ above. 
\par
In particular we see that $L(s)$ has a pole of order $2$ at $s=1$
(by [31] since $\lambda\not=1/4$ so that $F$ is a cusp form on
$GL(3)$). 
We can now proceed along classical lines: let $U>1$ (to be chosen
later) and let $\psi$ be a $C^{\infty}$ function on $[0,+\infty($ such
that $0\leq \psi\leq 1$ and 
$$
\psi(x)=\cases{
1& if  $0\leq x\leq 1-U^{-1}$,\cr
0& if $x\geq 1+U^{-1}$.
}
$$
The Mellin transform $\hat{\psi}(s)$ is holomorphic for $\sigma>0$, it
satisfies $\hat{\psi}(s)=s^{-1}+O(|\sigma|U^{-1})$ and by integration
by parts
$$
\hat{\psi}(s)\ll U^A(1+|t|)^{-A-1}
\leqno{(3.6)}
$$
for $\sigma\geq 1/2$ and for any $A>0$, the implied constant depending
on $A$ and $\psi$ only. 
\par
For suitable choices (say $\psi_+$ and $\psi_-$) of $\psi$ we get
$$
\sumb_{n\geq 1}{|\lambda_f(n)|^4\psi_-(n/x)}
\leq \sumb_{n\leq x}{|\lambda_f(n)|^4}\leq 
\sumb_{n\geq 1}{|\lambda_f(n)|^4\psi_+(n/x)}.
$$
Thus it is enough to prove (3.4) for a sum weighted by
$\psi(n/x)$. We have
$$
\eqalign{
\sumb_{n\geq 1}{|\lambda_f(n)|^4\psi(n/x)}&=
\frac{1}{2\pi i}\int_{(3)}{L(s)x^s\hat{\psi}(s)ds}\cr 
&=\frac{1}{2\pi i}\int_{(3)}{\zeta(s)L(F,s)^2L(F\otimes
F,s)H(s)H_1(s)\hat{\psi}(s)x^sds}.}
$$
For any fixed $\alpha>\frac{23}{32}$ we can move the line of
integration (the three $L$-functions are polynomially bounded in
vertical strips and $\hat{\psi}$ decays rapidly) to
$\sigma=\alpha$. We 
pass the double pole at $s=1$ with residue of the form
$$
d_fx\log x+e_fx+O(x(\log x)(\Lambda N)^{\eps}U^{-1})
$$
with $d_f=L(F,1)H(1)H_1(1)\hbox{res}_{s=1}L(F\otimes
F,s)>0$, $d_f$ and $e_f$ being estimated by [32] to get
$$
d_f\ll (\Lambda N)^{\eps},\quad e_f\ll (\Lambda N)^{\eps}
$$
for any $\eps>0$, the implied constant depending only on $\eps$.
\par
Now the integral on $\sigma=\alpha<1$ is estimated using 
$H(s)H_1(s)\ll 1$, the uniform convexity bound for automorphic
$L$-functions (see e.g. [24, \S 5.12]) yielding 
$$
L(s)\ll 
(1+|t|)^{8(1-\alpha)+\eps}\Lambda^{5(1-\alpha)+\eps}N^{8(1-\alpha)+\eps} 
$$
for the product of the three $L$-functions, the implied constant
depends only on $\alpha$ and $\eps$. Then
(3.6) with 
$A=8(1-\alpha)+1/2$ (to get an absolutely convergent integral) yields 
$$
\eqalign{
\sumb_{n\geq 1}{|\lambda_f(n)|^4\psi(n/x)}&=
d_fx\log x+e_fx+O(x(\log x)(\Lambda N)^{\eps}U^{-1})
+\cr
&\quad \quad O(x^{\alpha}\Lambda^{5(1-\alpha)+\eps}N^{8(1-\alpha)+\eps}
U^{8(1-\alpha)+1/2+\eps}).}
$$
Without trying to optimize, we take $U$ so that
$xU^{-1}=x^{\alpha}U^{8(1-\alpha)+1/2}$, which gives
$$
\sumb_{n\geq 1}{|\lambda_f(n)|^4\psi(n/x)}=
d_fx\log x+e_fx+O(\Lambda^{5(1-\alpha)+\eps}N^{8(1-\alpha)+\eps}
x^{\beta+\eps})
$$
with
$$
\beta=1-\frac{2(1-\alpha)}{16(1-\alpha)+3}.
$$
Taking $\alpha=\frac{23}{32}+\eps$, we get (3.4), up to renaming
$\eps$.
\par
For holomorphic forms, we proceed in the a slightly different manner.
First since we have a nebentypus we use the adjoint square instead of
the symmetric square in proving the analogue of (3.3), namely
$$
\sumb_{n\leq x}{|\lambda_f(n)|^2}=c_f x^k +
O(x^{k-1/5+\eps}(kN)^{1/2+\eps})
\leqno{(3.7)}
$$
with $c_f\gg (\log kN)^{-1}$ (by Goldfeld, Hoffstein and Lieman, see the
Appendix to [22]). Secondly we can avoid proving the analogue of
(3.4), for which we require only an upper bound, by means of
the Ramanujan-Petersson bound (proved by Deligne)
$$
|\lambda_f(n)|\leq d(n)n^{(k-1)/2},
$$
where $d(n)$ is the divisor function. In fact it is more efficient
then to argue with the third power moment, and use H\"older's
Inequality with $(p,q)=(3,2/3)$ for the final estimates:
$$
\sumb_{x<n\leq x+y}{|\lambda_f(n)|^2}
\leq \eta |S| + \sumb_{x\in L}{|\lambda_f(n)|^2}
\leq \eta y + |L|^{1/3} \Bigl(
\sumb_{x<n\leq x+y}{|\lambda_f(n)|^3}\Bigr)^{2/3}.
\leqno{(3.8)}
$$
We have
$$
\sumb_{x<n\leq x+y}{|\lambda_f(n)|^3}\leq 
\sum_{x<n\leq x+y}{d(n)^3n^{3(k-1)/2}},
$$
and estimating this is classical. Here are the main steps
for completeness. 
The generating Dirichlet series for $d(n)^3$ is
$$
L_1(s):=\sum_{n\geq 1}{d(n)^3n^{-s}}=
\zeta(s)^4\prod_p{(1+4p^{-s}+p^{-2s})}
=\zeta(s)^{8}H_2(s)
$$
where $H_2$ is absolutely convergent, hence holomorphic, for
$\sigma>\frac{1}{2}$. Say it has coefficients $\alpha(n)$, and
$\zeta(s)^8$ has coefficients $d_8(n)$. By [21, Th. 2] we have
$$
\sum_{n\leq x}{d_8(n)}=xP(\log x)+O(x^{5/8+\eps}) 
$$
where $P$ is some polynomial of degree $7$. Hence
$$
\eqalign{
\sum_{n\leq x}{d(n)^3}&=\sum_{b\leq x}{\alpha(b)
\sum_{a\leq x/b}{d_8(a)}}\cr
&=xP_1(\log x) +O(x^{5/8+\eps})
}
$$
for some polynomial $P_1$ of degree $7$ since $H_2(5/8)$
is absolutely convergent.
By partial summation we get
$$
\sum_{n\leq x}{d(n)^3n^{3(k-1)/2}}=x^{3(k-1)/2+1}P_2(\log x)+
O(x^{3(k-1)/2+5/8+\eps}),
$$
hence the result follows using (3.7), (3.8) since $5/8<4/5$.
\hfill$\square$

\smallskip

{\sl Remark 2.} We see that this method provides $n$ where a lower
bound for $\lambda_f(n)$ holds, and this also seems very hard to get
by purely algebraic techniques. In applications to analytic number theory,
this can be of crucial importance; see for instance [8], [9]. In these
papers the question is somewhat different: one needs to find very small $n$,
compared to some large parameter $x$ (say $n\ll x^{\eps}$), such that
$\lambda_f(n)$ is not too small, and this is solved by using the trick
of Iwaniec that for any prime $p\,\nmid N$, we have
$\lambda_f(p)^2-\lambda_f(p^2)=1$, so one of $\lambda_f(p)$,
$\lambda_f(p^2)$ is at least $1/\sqrt{2}$ in absolute value, and $p^2$
remains small enough for the application in mind.
\medskip

There is a strong contrast between the proof of Proposition 3, which
depends on quite deep 
analytic properties of $L$-functions, and the algebraic approach of
the previous section, where not even convergence mattered! It
is clear that one can extend Proposition 2 to any cuspidal
automorphic form on $GL(n)/\Q$ using its Rankin-Selberg convolution
(compare [9]), but Proposition 3 requires either that $f$ satisfies
the Ramanujan-Petersson conjecture, or that the adjoint square
be automorphic (in which case there is also is a bound of the type
$|\alpha_p|\leq p^{\theta}$ with $\theta<\frac{1}{4}$ for the local
parameters of $f$ at unramified primes). This is not known for $n\geq
3$. 
\par
\medskip
It is natural to ask if the property in Lemma 2.1 holds for primitive
Maass forms. If the eigenvalue is $\lambda=1/4$, conjecturally the Fourier
coefficients still generate a number field, and in this case the proof
goes through without change. If $\lambda\not=1/4$, the field
$K_f=\Q(\lambda_f(n),\chi(n))$ is not expected to be a number
field. However we still see that if Lemma 2.1 is false for $f$, then  
$\Q(\alpha_p,\beta_p,\chi(p))\cap \Q^{ab}$ is an infinite extension of
$\Q$, where  
$\Q^{ab}$ is the cyclotomic field generated by all roots of unity.
This does not sound very likely, as the field generated by the local
roots $\alpha_p$, $\beta_p$ could be expected to be
mostly transcendental, but it is certainly beyond proof or
disproof today! (The corresponding fact is true however, for the field
$K_t=\Q(2^{it},3^{it},\ldots,p^{it},\ldots)$ generated by the local
roots of the Eisenstein series $E(z,\dm+it)$, for $SL(2,\Z)$
say, for all $t\in \R$ except maybe those in a countable set;
it doesn't seem easy to decide if the latter is really empty, but this
would follow from Schanuel's Conjecture, as observed by B. Poonen).
\par
One is tempted to confront this with the famous ``optimistic''
question of Katz ([25, p.15]): is
$$
L(s)=\prod_{p}{(1-S(1,1;p)p^{-s}+p^{1-2s})^{-1}}=\sum_{n\geq 1}
{\lambda_S(n)n^{-s}} 
$$
the $L$-function of a (primitive) Maass form (of weight $2$), even up
to finitely many factors, where $S(1,1;p)$ denotes the usual
Kloosterman sums? Note 
that $S(1,1;p)$ generates the maximal real 
subfield of the field of $p$-th roots of unity, so in
this case the field generated by $\lambda_S(p)$ is an infinite
algebraic extension of $\Q$. However we can prove
the analogue of Lehmer's conjecture for this Dirichlet series! (Of
course, the answer to Katz's question is widely expected to be ``No'',
see [5] for some strong evidence).

\proclaim Proposition 4. 
For any $n\geq 1$, we have $\lambda_S(n)\not=0$.

\noindent{\sl Proof}. 
We give two proofs (suggested by Katz and simpler than our original
argument). We need to show that $\lambda_S(p^{\nu})\not=0$ for $p$ prime
and $\nu\geq 0$. For the first argument, consider the Euler factor at
$p$ as a rational function of $X=p^{-s}$ with coefficients in the
cyclotomic field $\Q(e(1/p))$. It is congruent (modulo the ideal
generated by $p$) to
$$
\frac{1}{1-S(1,1;p)X}=\sum_{\nu}{S(1,1;p)^{\nu}X^{\nu}}.
$$
Thus the result follows from the well-known fact that $S(1,1;p)$ is
non-zero modulo $p$, in fact we have
$$
S(1,1;p)\equiv -1\pmod{1-e(1/p)},
$$
and the prime ideal $1-e(1/p)$ divides $p$.
\par
For the other argument, notice that since the form of the Euler
product is the same as for a 
holomorphic form of weight $2$, we must show that
$\alpha_p/\beta_p$ is not a root of unity, where $\alpha_p$ and
$\beta_p$ satisfy
$$
\alpha_p+\beta_p=S(1,1;p)\quad \hbox{ and } \quad\alpha_p\beta_p=p.
$$
Hence the product $\alpha_b\beta_p$ is divisible by $1-e(1/p)$,
whereas by the congruence above, the sum is invertible modulo
$1-e(1/p)$. This means one of 
$\alpha_p$, $\beta_p$ must also be invertible while the other is not,
which implies that the ratio $\alpha_p/\beta_p$ is not a $p$-unit,
hence not a root of unity. 
\hfill$\square$
\smallskip

It is probably possible to derive a fancy proof of this proposition
(more amenable to generalizations, if desired)
using ideas as in [14], Lemma 4.9, applied to some Kloosterman/Gauss
sum sheaves on $\G_m/\F_p$ with traces of Frobenius at $\alpha\in
\G_m(\F_p)$ given by both sides of (3.9). Note also that if $\nu\geq
1$ and $p$ is odd we do have (see e.g. [23, Lemma 4.1])
$$
S(1,1;p^{2\nu})=p^{\nu}\Bigl(e\Bigl(\frac{2}{p^{2\nu}}\Bigr)
+e\Bigl(\frac{-2}{p^{2\nu}}\Bigr)\Bigr),
$$
so Proposition 4 is special to Kloosterman sums with prime modulus.

\smallskip
\vskip 5mm

\noindent{\bf \S\ 4. Applications of $\ciB$-free numbers}
\medskip

We now come to the technical heart of this paper where we consider the
original question of proving (1.5) for a cusp form $f\in
S_k(N,\chi)$, not in the space spanned by CM forms. Recall that Balog
and Ono [2] proved (1.5) for  
$y=x^{17/41+\eps}$, $\eps>0$ being arbitrary. 
It is interesting to look for smaller exponents, in particular since
it is natural to expect that $y=x^{\eps}$ should be sufficient. (By a
result of Plaksin [35] on $\ciB$-free numbers, this is true for almost
all $n$). For one very natural $f$, namely the Ramanujan $\Delta$
function with coefficients $\tau(n)$, a famous conjecture of Lehmer
[29] says that $\tau(n)\not=0$ for any $n\geq 1$. 
\par
Since this problem seems very difficult, approaching it by means of
conditional statements based on solid conjectures is also desirable.
Very recently Alkan [1] gave two such results:
he showed that the exponent $17/41$ can be reduced to $69/169$
and $1/3$ ([1], Theorems 3 and 4)
under the generalised Riemann hypothesis (GRH) for Dedekind zeta-function and
the Lang-Trotter conjecture [28], respectively.
\goodbreak
\smallskip
We will prove a number of results improving the previously known
statements, both conditional and unconditional. The following is a
general bound, where we recall that 
$\ciP_{f,1}$ is defined in (2.6):

\proclaim Theorem 1.
Suppose that $k\ge 2$ and $f\in S_k^*(N, \chi)$ is
a primitive form not of CM type 
such that
$$\big|\ciP_{f, 1}\cap [1, x]\big|
\ll_f x^\rho {(\log\log x)^{\Psi_\rho}\over (\log x)^{\Theta_\rho}}
\qquad(x\ge 2),
\leqno(4.1)$$
where $\rho\in [0, 1]$ and $\Theta_\rho, \Psi_\rho$ are real
constants such that $\Theta_1>1$.
Define
$$\theta(\rho)
= \cases{
{1\over 4}               & if $0\leq \rho\le {1\over 3}$,
\cr\noalign{\medskip}
{10\rho\over 19\rho + 7}    & if ${1\over 3}<\rho\le
{9\over 17}$,
\cr\noalign{\medskip}
{3\rho\over 4\rho+3}         & if ${9\over 17}<\rho\le
{15\over 28}$,
\cr\noalign{\medskip}
{5\over 16}                 & if ${15\over 28}<\rho\le
{5\over 8}$,
\cr\noalign{\medskip}
{22\rho\over 24\rho + 29}    & if ${5\over 8}<\rho\le
{9\over 10}$,
\cr\noalign{\medskip}
{7\rho\over 9\rho + 8}       & if ${9\over 10}<\rho\le 1$,
\cr}$$
For every $\varepsilon>0$, $x\ge x_0(f, \varepsilon)$ and $y\ge 
x^{\theta(\rho)+\eps}$, we have
$$|\{n\,\mid\, x<n\le x+y\hbox{ and }\lambda_f(n)\not=0\}|\gg_{f,
\varepsilon} y. 
\leqno(4.2)$$
In particular for any $\varepsilon>0$ and all $n\ge 1$, we have
$$i_f(n)\ll_{f, \varepsilon} n^{\theta(\rho)+\eps}.
\leqno(4.3)$$

\smallskip

Theorem 1 follows immediately by multiplicativity from Corollary 10
below which 
gives a more effective treatment for
$\ciB$-free numbers in short intervals, applied with
$$
\ciP=\{p\,\mid\, p\,\mid N\hbox{ or }
\lambda_f(p)=0\}.
$$
The new ideas and new ingredients will be explained in \S\ 5.
\smallskip

According to (1.2),
the hypothesis (4.1) holds with $(\rho, \Theta_\rho, \Psi_\rho) = (1,
1+\delta, 0)$ for any $\delta<\dm$.
Thus, applying this result and Lemma 2.4, we immediately obtain an
improvement of the result of Balog and Ono. 

\proclaim Corollary 1.
Suppose that $k\ge 2$ and $f\in S_k(N, \chi)$ is
not in the space spanned by CM forms.
Then for any $\varepsilon>0$, $x\ge x_0(f, \varepsilon)$ and $y\ge
x^{7/17+\varepsilon}$,
we have
$$|\{n\,\mid\, x<n\le x+y\hbox{ and }\lambda_f(n)\not=0\}|\gg_{f,
\varepsilon} y.$$ 
In particular
$$i_f(n)\ll_{f, \varepsilon} n^{7/17+\varepsilon}.$$

In proving this we do not exploit Lemma 2.1 (so we could claim that
we obtain the correct proportion of squarefree numbers if $f$ is
primitive). It can be used to simplify the proof, as we'll see, but it
does not influence the strength of the exponent. This is mainly due to
the fact that we have $\rho=1$, and when $\rho$ is close to $1$ we do not
succeed in getting better results by not imposing the numbers to be
squarefree.  
\par
However, if one can get $\rho$ quite small, e.g. smaller than the
current best results about squarefree numbers in short intervals (see
[13]), it is clear that using Lemma 2.1 will yield an improvement.
So consider the set of prime numbers
$$\ciP_f^* := \{p\,\mid\,p|N \}\cup
\mathop\cup_{\nu=1}^\infty \ciP_{f, \nu},$$
where as before
$$\ciP_{f, \nu}
= \{p\,\mid\,p\,\nmid N \,\, {\rm and} \,\, \lambda_f(p^\nu)=0\}.$$
Clearly $\lambda_f(n)\not=0$ (for $(n,N)=1$) if and
only if $n$ is 
$\ciP_f^*$-free. 
We then have the following result:

\proclaim Theorem 2.
Assume that $k\ge 2$, $f\in S_k^*(N, \chi)$ is
not a CM form, and that
$$\big|\ciP_f^*\cap [1, x]\big|
\ll_f x^\rho {(\log\log x)^{\Psi_\rho}\over (\log x)^{\Theta_\rho}}
\qquad(x\ge 2),
\leqno(4.4)$$
where $\rho\in [0, 1]$ and $\Theta_\rho, \Psi_\rho$ are real
constants such that $\Theta_1>1$.
Then the inequalities (4.2) and (4.3) hold with $\theta(\rho) = 
\rho/(1+\rho)$ for $0\leq \rho\le 1$.

This theorem gives a better exponent than Theorem 1 when $\rho\le
{1\over 3}$ under a slightly stronger hypothesis than (4.1). However recall
from Lemma 2.2 that the hypotheses (4.1) and (4.4) are in fact equivalent
when $k$ is even. It is of course particularly interesting that this
new exponent tends towards $0$ when $\rho\to 0$. As for Theorem 1,
this result follows directly by multiplicativity from the
corresponding result for $\ciB$-free numbers, Proposition 9 below,
where this time $\ciP=\ciP_f^*$.
\par
\smallskip
Another consequence of Lemma 2.3 and Corollary 10 is an extension to
all symmetric powers:

\proclaim Corollary 2. 
Let $k\geq 2$ and $f\in S_k^*(N,\chi)$ which is not a CM form. Let
$m\geq 1$ and define the unramified $m$-th symmetric power
$L$-function of $f$ by
$$
L_{nr}(\symm{m} f,s)=\prod_{p\nmid N}{
\prod_{0\leq j\leq m}{(1-\alpha_p^j\beta_p^{m-j}p^{-s})^{-1}}}
=\sum_{n\geq 1}{\lambda_f^{(m)}(n)n^{-s}}.
$$
Then for any $\varepsilon>0$, $x\ge x_0(f, \varepsilon)$ and $y\ge
x^{7/17+\varepsilon}$,
we have
$$
|\{n\,\mid\, x<n\le x+y\hbox{ and }\lambda_f^{(m)}(n)\not=0\}|\gg_{f,m,
\varepsilon} y,
$$
and in particular $i_{{\rm Sym}^m f}(n)\ll_{f, \varepsilon, m}
n^{7/17+\varepsilon}$ for $n\geq 1$.

\noindent{\sl Proof.} For $p\nmid N$ prime, we have 
$\lambda_f^{(m)}(p)=\lambda_f(p^{m})$. Hence by Lemma 2.3 we derive
$$
|\{p\leq x\,\mid\, p\nmid N\hbox{ and } \lambda_f^{(m)}(p)=0\}|
\ll x(\log x)^{-1-\delta}
$$
for any $\delta<\dm$. By multiplicativity and Corollary 10 below, the result
follows. \hfill$\square$

\smallskip
Note we do not need the automorphy of $\symm{m} f$ (which is known only
for $m\leq 4$).
\par
\smallskip
The hypothesis (4.1) is known only with $\rho=1$, with the one
exception of primitive forms $f\in S_2^*(N)$ with integral
coefficients. Those are associated to elliptic curves over $\Q$, and
Elkies [10] has proved that (4.1) (or (4.4)) holds with $\rho=3/4$,
$\Theta=\Psi=0$. Theorem 1 is still better for this value of $\rho$ than
Theorem 2 and we get:
\smallskip

\proclaim Corollary 3.
Let $E/\Q$ be an elliptic curve without complex multiplication and let
$f$ be the associated primitive form. Then
for every $\varepsilon>0$, $x\ge x_0(E, \varepsilon)$ and $y\ge
x^{33/94+\varepsilon}$,
we have
$$|\{n\,\mid\, x<n\le x+y\hbox{ and }
\lambda_f(n)\not=0\}|\gg_{E, \varepsilon} y.$$
In particular for any $\varepsilon>0$ and all $n\ge 1$, we have
$$i_{f}(n)\ll_{E, \varepsilon} n^{33/94+\varepsilon}.$$

This improves Theorem 2 of [1], which requires $69/169$ in
place of $33/94$.
\par
\smallskip
Some well-known conjectures imply that (4.1) holds for smaller values
of $\rho$. For example, Serre ([43, $(182)_{\rm R}$]) showed that the GRH for 
Dedekind zeta-functions
implies (4.1) with $(\rho, \Theta_\rho, \Psi_\rho) = ({3\over 4}, 0, 0)$.
Lang and Trotter [28] formulated a conjecture for the size of the set 
$\ciP_{f, 1}$,
in the case where $f$ is associated to an elliptic curve over $\Q$.
This, if true, implies for these forms an estimate (4.1) with $(\rho,
\Theta_\rho, \Psi_\rho) = ({1\over 2}, 1, 0)$.
Generalizations of the Lang-Trotter conjecture (see e.g. Murty's
version [33], especially Conjecture 3.4) imply that
if $k\ge 2$ and
$f\in S_k^*(N, \chi)$ is not of CM type, then
we have (4.1) with
$$(\rho, \Theta_\rho, \Psi_\rho) = \cases{
(\dm, 1, 0)              & if $k=2$ and $[F_f : \Q] = 2$,
\cr\noalign{\smallskip}
(0, 0, 1)                & if $k=2$ and $[F_f : \Q] = 3$
\cr\noalign{\smallskip}
& or $k=3$ and $[F_f : \Q] = 2$,
\cr\noalign{\smallskip}
(0, 0, 0)                & otherwise,
\cr}
\leqno(4.5)$$
where $F_f$ is the stable trace field (see \S\ 2 and \S\ 3 of
[33]).
\goodbreak
\smallskip
Applying Theorem 1, we get the following conditional result,
which improves Theorem 1 of [1].

\proclaim Corollary 4.
Suppose that $k\ge 2$ and $f\in S_k(N, \chi)$ is
not in the space spanned by CM forms.

\vskip -2mm
{\sl {\rm (i)}
Under the GRH for Dedekind zeta-function,
the exponent $7/17$ of Corollary 1 can be further improved to $33/94$.

{\rm (ii)}
Under the generalized Lang-Trotter conjecture,
the exponent $7/17$ can be further improved to $10/33$ if
$k=[F_f : \Q] = 2$,
and to $1/4$ otherwise.}

\medskip
We can apply Theorem 2 instead if $f$ satisfies the assumptions of
Lemma 2.2, but it is just as simple to extend the Lang-Trotter type
conjectures to deal with the sets $\ciP_{f,\nu}$ for any $\nu\geq
1$. The heuristics which lead to these conjectures, based on Deligne's
estimate $|\alpha_p|=|\beta_p|=p^{(k-1)/2}$, suggest the following:

\proclaim Conjecture 1.
Let $\nu\ge 1$ be any integer.
If $k\ge 2$ and $f\in S_k^*(N,\chi)$ is not of CM type,
then
$$\big|\ciP_{f, \nu}\cap [1, x]\big|
\ll_f x^\rho {(\log\log x)^{\Psi_\rho}\over (\log x)^{\Theta_\rho}}
\qquad(x\ge 2)$$
with $(\rho, \Theta_\rho, \Psi_\rho)=(\dm,1,0)$ if $k=2$, $(0,0,1)$ if
$k=3$ and $(0,0,0)$ if $k\geq 4$.

We only state upper bounds, but one could
propose a more precise statement, which involves looking at the
possibility of $f$ having ``extra twists'' 
and eliminating the all but finitely many $\nu$ for
which $\ciP_{f, \nu}$ is empty.
About this conjecture, recall that even under GRH, one can not get a
better general result towards the Lang-Trotter conjecture than
$$\big|\ciP_{f, 1}\cap [1, x]\big|
\ll_f x^{3/4}$$
for $f$ of weight $k\geq 2$. The exponent is the same for all weights, so
this gets worse (compared to what we expect) as $k$ grows. In
particular, this conjecture for $k\geq 3$ seems hopeless for the time
being. Lemma 2.1 implies:

\proclaim Corollary 5.
Let $k\ge 2$ and $f\in S_k^*(N,\chi)$ not of CM type.
Assuming Conjecture 1 for $f$, the inequality (4.4) holds with
$(\rho, \Theta_\rho, \Psi_\rho)$ given by $(\rho, \Theta_\rho,
\Psi_\rho)=(\dm,1,0)$ if $k=2$, $(0,0,1)$ if 
$k=3$ and $(0,0,0)$ if $k\geq 4$.

\smallskip
As applications (or cautionary tale...), here 
are some very impressive-looking results.

\proclaim Corollary 6.
Suppose that $k\ge 3$ and $f\in S_k(N, \chi)$ is
not in the space spanned by CM type. If Conjecture 1 holds for all
primitive forms, then the exponent $7/17$ of Corollary 1 can be
improved to $0$. If $k\geq 4$ and $f$ is primitive, then there exists
$M\geq 1$ such that $(n,M)=1$ implies $\lambda_f(n)\not=0$.

\smallskip
Specializing to the Ramanujan $\tau$-function, which is integer
valued, Lemma 2.2 allows us to deduce the following result (implicit
in [43]):

\proclaim Corollary 7.
Assume Conjecture 1, or equivalently the generalized Lang-Trotter
conjecture, for $f=\Delta\in S^*_{12}(1)$.
There exists $P\geq 1$ such that
$\tau(n)=0$ if and only if $(n,P^{\infty})$ is a square, i.e. if and
only if $v_p(n)$ is even for $p\mid P$.
In particular $i_{\Delta}(n)\leq P$ for $n\ge 1$, and
for all $x\geq 2$ and $y\geq 1$, we have
$$
|\{n\,\mid\,x<n\leq x+y\,\,\hbox{and} \,\, \tau(n)\not=0\}|
= \prod_{p\mid P}{\Bigl(1+\frac{1}{p}\Bigr)^{-1}} y +O((\log
(x+y))^{\omega(P)}) \geq \frac{\varphi(P)}{P}y+O(1)
$$
where the implied constant is absolute and $\omega(P)$ is the number
of prime divisors of $P$.
\medskip

\noindent{\sl Proof}.
The first statement is the rephrasing of Lemma 2.2 and Conjecture 1 in
this case. Notice that $\tau(n)\not=0$ if $(n,P)=1$ so
$i_{\Delta}(n)\leq P$ follows (an interval of length $P$ contains
elements prime to $P$) as does the last inequality by trivial
counting. For the asymptotic, write
$$
\eqalign{
|\{n\leq x\,\mid\, \tau(n)\not=0\}|&=
\sum_{\stacksum{d\mid P^{\infty}}{d\leq x}}{\lambda(d)\sum_{n\leq x/d}{1}}
}
$$
where $\lambda(n)$ is the Liouville function,
i.e. $\lambda(p^k)=(-1)^k$. Since 
$$
\sum_{d\mid P^{\infty}}{\frac{\lambda(d)}{d}}=
\prod_{p\mid P}{(1+p^{-1})^{-1}},
$$
we get the result after elementary estimates.
\hfill$\square$

\medskip
This is of course trivial and of little practical significance towards the
Lehmer conjecture.

\smallskip
\vskip 5mm
\goodbreak
\noindent{\bf \S\ 5. Multiple exponential sums and bilinear forms}
\medskip

This section is devoted to the study of multiple exponential sums and
bilinear forms, which will be used in the proofs of our results on
$\ciB$-free numbers in the next sections, but are also of independent
interest. 
We begin by investigating a double exponential sum of type II:
$$S(M, N) := \sum_{m\sim M}\sum_{n\sim N} \varphi_m \psi_n
e\bigg(X{m^\alpha n^\beta\over M^\alpha N^\beta}\bigg),$$
where $e(t) := \exp\{2\pi it\}$,
$X>0$,
$M\ge 1$,
$N\ge 1$,
$|\varphi_m|\le 1$,
$|\psi_n|\le 1$,
$\alpha, \beta\in \R$
and $m\sim M$ means $M\le m<2M$.
Such a sum occurs in many arithmetic problems
and is studied by many authors (for example, [15] and [40]).
We shall estimate this sum by the method of Fouvry \& Iwaniec [15]
together with the refinement of Robert \& Sargos [39].
When $X<N^2$, we need to use an idea in [40].

The following result is an improvement of Theorem 4 in [15] and 
Theorem 10 in [40].

\proclaim Proposition 5.
If $\alpha, \beta\in \R\sset\{0,1\}$, then for any $\varepsilon>0$ we have
$$S(M, N)
\ll \big\{
   (XM^6N^6)^{1/8}
+ M^{1/2}N
+ MN^{3/4}
+ X^{-1/2} MN
\big\}(MN)^\varepsilon.$$

\noindent{\sl Proof}.
We shall distinguish two cases.

\smallskip

{A. \sl The case of $X\ge N^2$}

\smallskip

By applying twice the Cauchy-Schwarz' inequality, it follows that
$$|S(M, N)|^4
\le (M N)^2 \sum_{n_1\sim N} \sum_{n_2\sim N}
\sum_{m_1\sim M} \sum_{m_2\sim M}
e\bigg(X{(m_1^\alpha-m_2^\alpha)(n_1^\beta-n_2^\beta)\over M^\alpha 
N^\beta}\bigg).$$

The double large sieve inequality ([15], Proposition 1) with the choice of
$${\cal X} = \big\{(m_1^\alpha-m_2^\alpha)/M^\alpha\big\}_{m_1, m_2\sim M}
\quad{\rm and}\quad
{\cal Y} = \big\{(n_1^\beta-n_2^\beta)/N^\beta\big\}_{n_1, n_2\sim N}$$
leads to the following estimate
$$|S(M, N)|^8
\ll X (MN)^4 {\cal N}(M, 1/X) {\cal N}(N, 1/X),
\leqno(5.1)$$
where ${\cal N}(M, \Delta)$ is the number of quadruplets
$(m_1, m_2, m_3, m_4)\in \{M+1, \dots, 2M\}^4$ satisfying
$$\big|m_1^\alpha+m_2^\alpha-m_3^\alpha-m_4^\alpha\big|\le \Delta M^\alpha.$$
According to Theorem 2 of [39], we have
$${\cal N}(M, 1/X)
\ll \big(M^2 + X^{-1} M^4\big) M^\varepsilon.$$
Inserting this into (5.1) and
simplifying the estimate obtained by using the hypothesis $X\ge N^2$,
we find that
$$S(M, N)
\ll \big\{
   (XM^6N^6)^{1/8}
+ MN^{3/4}
\big\}(MN)^\varepsilon.$$

\smallskip

{B. \sl The case of $X\le N^2$}

\smallskip

By Lemma 2.1 of [40], we deduce that, for any $Q\in [1, M^{1-\varepsilon}]$,
$$|S(M, N)|^2\ll (M N)^2 Q^{-1} + M N Q^{-1} (\log M) \max_{1\le
Q_1\le Q} |S(Q_1)|,
\leqno(5.2)$$
where
$$S(Q_1):=
\sum_{q\sim Q_1}
\sum_{m\sim M} \varphi_{m, q}
\sum_{n\sim N}
e\bigg(X'{t(m, q) n^\beta\over T N^\beta}\bigg)
\leqno(5.3)$$
and
$$t(m, q) := (m+q)^\alpha - m^\alpha,
\qquad
T := M^{\alpha-1}Q_1,
\qquad
X' := XM^{-1}Q_1.$$

If $N':= X'/N\ge \dm$,
applying Lemma 2.2 of [40] to the sum over $n$ yields
$$\sum_{n\sim N} \!\!
e\bigg(X'{t(m, q) n^\beta\over T N^\beta}\bigg)
\ll X'^{-1/2} N
\!\!\!\! \sum_{n'\in I(m, q)} \!\!\!
w_{n'} e\bigg(\tilde \beta X'{u(m, q) {n'}^{\beta_1}\over U
N'^{\beta_1} }\bigg)
+ R_1 + R_2 + \log N,
\leqno(5.4)$$
where
$$\eqalign{
I(m, q)
& := \big[c_1 X t(m, q) M^{-\alpha} N^{-1}, \,
c_2 X t(m, q) M^{-\alpha} N^{-1}\big],
\cr\noalign{\smallskip}
R_j
& :=
\min\big\{X'^{-1/2} N, \, 1/ \|c'_j X M^{-\alpha} N^{-1} t(m, q)\|\big\},
\cr}$$
$u(m, q) := t(m, q)^{1/(1-\beta)}$,
$U := T^{1/(1-\beta)}$,
$\beta_1 := \beta/(\beta-1)$,
$\tilde \beta := |1-\beta| |\beta|^{-\beta_1}$,
$|w_{n'}|\le 1$,
and $c_j = c_j(\beta)$, $c'_j = c'_j(\beta)$ are some suitable constants.
Inserting into (5.3), using Lemma 2.5 of [40] to eliminate
multiplicative restrictions
and using Lemma 2.3 of [40] with $n=m$ to estimate the related
error terms, we find
$$S(Q_1)
\ll X'^{-1/2} N \int_{-\infty}^{+\infty} \Xi(r) S(Q_1,r) \d r
+ \big\{( X M^{-1} Q_1^3 )^{1/2}
+ M Q_1\big\} (MN)^\varepsilon,$$
where $\Xi(r) := \max\{M, \, (\pi r)^{-1}, \, (\pi r)^{-2}\}$,
$\psi_{n'}(r) := w_{n'} e(r n')$ and
$$S(Q_1,r)
:= \sum_{q\sim Q_1} \sum_{m\sim M}
\bigg|\sum_{n'\sim N'} \psi_{n'}(r) \,
e\bigg(\tilde \beta X'{u(m, q) {n'}^{\beta_1}\over U 
N'^{\beta_1}}\bigg)\bigg|.$$

\goodbreak

If $X'/N\le \dm$, the Kusmin-Landau inequality (see e.g. [20], Theorem
2.1) implies 
$$S(Q_1)\ll X'^{-1} MNQ_1.$$
Thus we always have
$$\leqalignno{S(Q_1)
& \ll X'^{-1/2} N \int_{-\infty}^{+\infty} \Xi(r) S(Q_1,r) \d r
& (5.5)
\cr
& \qquad
+ \big\{(X M^{-1} Q_1^3 )^{1/2}
+ M Q_1
+ X'^{-1} MNQ_1\big\} (MN)^\varepsilon,
\cr}$$

Now by applying Cauchy-Schwarz' inequality, it follows that
$$|S(Q_1,r)|^2
\le MQ_1 \sum_{q\sim Q_1} \sum_{m\sim M}
\sum_{n'_1\sim N'} \sum_{n'_2\sim N'}
\psi_{n'_1}(r) \overline{\psi_{n'_2}(r)}
e\bigg(\tilde \beta X'{u(m, q)
\big({n'_1}^{\beta_1}-{n'_2}^{\beta_1}\big)\over U N'^{\beta_1}
}\bigg).$$
The double large sieve inequality with the choice of
$${\cal X} = \big\{u(m, q)/U\big\}_{m\sim M, q\sim Q_1}
\quad{\rm and}\quad
{\cal Y} =
\big\{({n'_1}^{\beta_1}-{n'_2}^{\beta_1})/{N'}^{\beta_1}\big\}_{n_1,
n_2\sim N'}$$
allows us to deduce
$$|S(Q_1,r)|^4
\ll (MQ_1)^2 X' {\cal N}^*(u, U; 1/X') {\cal N}(N'; 1/X')
\leqno(5.6)$$
uniformly for $r\in \R$,
where ${\cal N}^*(u, U; \Delta)$ is the number of quadruplets
$(m_1+q_1, m_2+q_2, m_1, m_2)$ such that $m_1, m_2\sim M$, $q_1,
q_2\sim Q_1$ and
$$\big|u(m_1, q_1) - u(m_2, q_2)\big|\le \Delta U,$$
and ${\cal N}(N; \Delta)$ is the number of quadruplets
$(n_1, n_2, n_3, n_4)\in \{N+1, \dots, 2N\}^4$ satisfying
$$\big|n_1^\beta + n_2^\beta - n_3^\beta - n_4^\beta \big|\le \Delta N^\beta.$$

Since
$$\big|u(m_1, q_1) - u(m_2, q_2)\big|
\asymp \big|t(m_1, q_1) - t(m_2, q_2)\big| T^{\beta/(1-\beta)},$$
we have, for some suitable constant $C>0$,
$${\cal N}^*(u, U; \Delta) = {\cal N}^*(t, T; C\Delta).$$
Noticing that
$$\big|t(m_1, q_1) - t(m_2, q_2)\big|
= \big|(m_1+q_1)^\alpha - m_1^\alpha - (m_2+q_2)^\alpha - m_2^\alpha\big|$$
and
$(m_1+q_1, m_1, m_2+q_2, m_2)\in \{M+1, \dots, 3M\}^4$, clearly we have
$${\cal N}^*(t, T; C\Delta)
\ll {\cal N}(M; C\Delta).$$
Thus Theorem 2 of [39] implies that
$$\eqalign{{\cal N}^*(u, U; 1/X')
& \ll \big(M^2+X'^{-1}M^4\big) M^\varepsilon,
\cr
{\cal N}(N', 1/X')
& \ll \big(N'^2 + X'^{-1} N'^4\big) N'^\varepsilon.
\cr}$$
Inserting these into (5.6), we obtain uniformly for $r\in \R$,
$$|S(Q_1,r)|^4
\ll \big\{
    X' M^4(N'Q_1)^2
+ (MN')^4Q_1^2
+ M^6(N'Q_1)^2
+ X'^{-1}M^6N'^4Q_1^2
\big\} (MN)^\varepsilon.$$
Combining this with (5.5), we find that
$$\eqalign{S(Q_1)
& \ll \big\{
    (X M^3 N^2 Q_1^3)^{1/4}
+ (X M Q_1^2)^{1/2}
+ (M^3 N Q_1)^{1/2}
\cr
& \qquad
+ (X M^5 Q_1^3)^{1/4}
+ ( X M^{-1} Q_1^3 )^{1/2}
+ M Q_1
+ X'^{-1} MNQ_1
\big\} (MN)^\varepsilon.
\cr}$$
Since $Q\le M^{1-\varepsilon}$, the fifth and sixth terms on the
right-hand side are superfluous.
Inserting the simplified estimate into (5.2) and taking $Q=
M^{1-\varepsilon}$, we find
$$\leqalignno{|S(M, N)|^2
& \ll \big\{
     MN^2
+ (X M^6 N^6)^{1/4}
+ (X M^3 N^2)^{1/2}
& (5.7)
\cr
& \qquad
+ (M^4 N^3)^{1/2}
+ (X M^8 N^4)^{1/4}
+ X^{-1} (MN)^2
\big\} (MN)^\varepsilon.
\cr}$$
Similarly by interchanging the role of $M$ and $N$, we also have
$$\leqalignno{|S(M, N)|^2
& \ll \big\{
     M^2N
+ (X M^6 N^6)^{1/4}
+ (X M^2 N^3)^{1/2}
& (5.8)
\cr
& \qquad
+ (M^3 N^4)^{1/2}
+ (X M^4 N^8)^{1/4}
+ X^{-1} (MN)^2
\big\} (MN)^\varepsilon.
\cr}$$
Now the required estimate follows from (5.7) if $X\le N^2$ and $M\le N$,
and from (5.8) when $X\le N^2$ and $M>N$.
This completes the proof.
\hfill
$\square$

\bigskip

Next as an application of Proposition 5, we consider a particular triple
exponential sum of type I:
$$S_I(H, M, N)
:= \sum_{h\sim H} \sum_{m\in I} \sum_{n\sim N}
\xi_h \psi_n
e\bigg(X{h^\beta m^{-\beta} n^\alpha \over H^\beta M^{-\beta}
N^\alpha}\bigg),$$
where $X>0$,
$H\ge 1$,
$M\ge 1$,
$N\ge 1$,
$|\xi_h|\le 1$,
$|\psi_n|\le 1$
and $I$ is a subinterval of $[M,2M]$.

\proclaim Corollary 8.
Let $\alpha, \beta\in \R$ satisfy
$\beta\not=-1, 0$ and $\alpha/(1+\beta)\not=0, 1$.
For any $\varepsilon>0$, we have
$$\leqalignno{S_I(H, M, N)
& \ll \big\{
    (X^3 H^6 M^2 N^6)^{1/8}
+ (X H^2 N)^{1/2}
+ HN
& (5.9)
\cr
& \quad
+ (XH^3M)^{1/4}N
+ X^{-1} HMN
\big\} (HMN)^\varepsilon,
\cr\noalign{\medskip}
S_I(H, M, N)
& \ll \big\{
    (X^{\kappa+\lambda} H^{1+\kappa+\lambda} M^{1+\kappa-\lambda}
N^{2+\kappa})^{1/{(2+2\kappa)}}
+ (XH^2N)^{1/2}
& (5.10)
\cr
& \quad
+ (HM)^{1/2}N
+ HN
+ X^{-1} HMN
\big\} (HMN)^\varepsilon,
\cr}$$
where $(\kappa, \lambda)$ is an exponent pair.

\noindent{\sl Proof}.
If $X/M\le \dm$, the Kusmin-Landau inequality implies
$$S_I(H, M, N)\ll X^{-1}HMN.$$

When $X/M>\dm$, applying Lemma 2.2 of [40] to the sum over $m$ and
using Lemma 2.3 of [40] with $n=n$ to estimate the related error
terms, we find
$$S_I(H, M, N)
\ll X^{-1/2} M S'
+ (HN
+ X^{1/2}H) \log M,$$
where
$$S'
:= \sum_{n\sim N} \sum_{h\sim H} \sum_{m'\in I'(h, n)}
\widetilde \psi_n \widetilde \xi_h \varphi_{m'}
e\bigg(\widetilde \alpha X
{h^{\beta'} m'^{\beta'} n^{\alpha'} \over H^{\beta'} M'^{\beta'}
N^{\alpha'}}\bigg),$$
where $I'(h, n)$ is a subinterval of $[M', 2M']$ with $M':= X / M$,
$\beta' := \beta/(1+\beta)$,
$\alpha' := \alpha/(1+\beta)$,
$\widetilde \alpha := |1+\beta| |\beta|^{\beta'}$,
$|\widetilde \xi_h|\le 1$,
$|\varphi_{m'}|\le 1$
and
$|\widetilde \psi_n|\le 1$.
Noticing that the exponents of $h$ and $m'$ are equal,
we can express this new triple sum as a double exponential sum over $(h', n)$
with $h' = h m'\in h I'(h, n)$.
We use Lemma 2.5 of [40] to relax the condition $h' = h m'\in h I'(h, n)$
to $h'\sim H' := H M' = X H / M$.
Finally applying Proposition 5 with $(M, N) = (N, H')$ yields the
desired estimate (5.9).
The last inequality follows from (3.11) of [30] with the choice of
$(H, M, N) = (H', 1, N)$.
This completes the proof.
\hfill
$\square$

\bigskip

Finally we study bilinear form of type I:
$$\sum_{m\sim M}\sum_{n\sim N} \psi_n \, r_{mn}(x,y),
\leqno(5.11)$$
where $|\psi_n|\le 1$ and
$$r_d(x,y)
:= \sum_{\scriptstyle x<n\le x+y\atop\scriptstyle d\,\mid\,n} 1 - {y\over d}.
\leqno(5.12)$$
In the sequel, $\varepsilon$ denotes an arbitrarily small positive number
and $\varepsilon'$ a constant multiple of $\varepsilon$,
which may be different in each occurrence.

\proclaim Corollary 9.
Let $y := x^\theta$ and $|\psi_n|\le 1$.
Then for any $\varepsilon>0$ we have
$$\sum_{m\sim M}\sum_{n\sim N} \psi_n \, r_{mn}(x,y)
\ll_\varepsilon yx^{-\varepsilon}
\leqno(5.13)$$
provided one of the following two conditions holds
$$\cases{
{1\over 3}<\theta\le {5\over 11},
\cr\noalign{\smallskip}
N\le y^{9/4}x^{-3/4-\varepsilon'},
\cr\noalign{\smallskip}
MN\le x^{1-\varepsilon'},
\cr}
\leqno(5.14)$$
or
$$\cases{
(\kappa+\lambda)/(1+2\kappa+2\lambda)<\theta\le
(\kappa+\lambda)/(2\kappa+\lambda),
\cr\noalign{\medskip}
N\le
y^{(1+2\kappa+2\lambda)/(1+\lambda)}x^{-(\kappa+\lambda)/
(1+\lambda)-\varepsilon'}, 
\cr\noalign{\vskip 1,3mm}
MN\le x^{1-\varepsilon'}.
\cr}
\leqno(5.15)$$

\noindent{\sl Proof}.
Without loss of generality, we can suppose that $MN\ge yx^{-\varepsilon}$.
By applying (5.9) of Corollary 8, we see that
$$\sum_{h\sim H} \sum_{m\sim M} \sum_{n\sim N}
\psi_n \, e\bigg({xh\over mn}\bigg)
\ll MNx^{-2\varepsilon},$$
provided
$$\textstyle {1\over 3}<\theta\le {5\over 11},
\qquad
H\le MNy^{-1}x^{3\varepsilon},
\qquad
N\le y^{9/4}x^{-3/4-\varepsilon'},
\qquad
MN\le x^{1-\varepsilon'}.$$
Combining this with Lemma 9 of [46] with the choice of $\varphi_m\equiv 1$,
we deduce (5.13) provided (5.14) holds.
The other one can be proved by using (5.10) of Corollary 8.
\hfill
$\square$

\medskip

A particular case of (5.11) -- linear forms (with $N=1$) --
will be needed in the proof of Corollary 10.

\proclaim Lemma 5.1.
Let $y := x^\theta$.
Then for any $\varepsilon>0$ we have
$$\sum_{m\sim M} r_m(x,y)
\ll_\varepsilon yx^{-\varepsilon}
\leqno(5.16)$$
provided one of the following two conditions holds
$$\leqalignno{
& \textstyle {1\over 4}<\theta\le {9\over 29}
\qquad\hbox{and}\qquad
M\le y^{19/7} x^{-3/7-\varepsilon'};
& (5.17)
\cr\noalign{\medskip}
& \textstyle {9\over 29}<\theta\le {1\over 2}
\qquad\hbox{and}\qquad
M\le y^{4/3} x^{-\varepsilon'}.
& (5.18)
\cr}$$

\noindent{\sl Proof}.
Without loss of generality, we can suppose that $M\ge yx^{-\varepsilon}$.
Theorem 1 of [38] allows us to write
$$\eqalign{\sum_{h\sim H} \bigg|\sum_{m\sim M} e\bigg({xh\over m}\bigg)\bigg|
& \ll \big\{(x^3H^{19}M^6)^{1/18}
+ (xH^6M)^{1/5}
+ H M^{3/4}
+ (x^{-1}H^2M^4)^{1/3}\big\} M^\varepsilon
\cr
& \ll Mx^{-2\varepsilon},
\cr}$$
provided one of the following two conditions holds
$$\textstyle {1\over 4}<\theta\le {9\over 29},
\qquad
H\le My^{-1}x^{3\varepsilon},
\qquad
M\le y^{19/7} x^{-3/7-\varepsilon'}$$
or
$$\textstyle {9\over 29}<\theta\le {1\over 2},
\qquad
H\le My^{-1}x^{3\varepsilon},
\qquad
M\le y^{4/3} x^{-\varepsilon'}.$$
This implies (5.16) if (5.17) or (5.18) holds.
\hfill
$\square$

\vskip 5mm

\noindent{\bf \S\ 6. $\ciB$-free numbers in short intervals}

\medskip

In this section we explain our new results about $\ciB$-free numbers. 
The notion of $\ciB$-free numbers, introduced by Erd\H os [11],
is a generalisation of square-free integers.
More precisely, let
$$\ciB = \{b_k\,\mid\,1<b_1<b_2<\cdots \,\,\}$$
be an infinite sequence of integers such that
$$\sum_{k=1}^\infty {1\over b_k} <\infty
\qquad{\rm and}\qquad
(b_j, b_k)=1
\quad
(j\not=k).
\leqno(6.1)$$
The $\ciB$-free numbers are the integers that are divisible by no
element of $\ciB$. We already mentioned that 
the existence of $\ciB$-free numbers in short intervals was 
proved by Erd\H os [11],
who showed that there is a constant $\theta\in (0, 1)$ such that the
short interval
$(x, x+x^\theta]$ with $x$ sufficiently large contains $\ciB$-free numbers.
Szemer\'edi [44] showed that $\theta = \dm+\varepsilon$ is admissible.
This result was further improved to
$$\eqalign{
& \theta = \textstyle {9\over 20}+\varepsilon
\qquad
\hbox{by Bantle \& Grupp [3]},
\cr
& \theta = \textstyle {5\over 12}+\varepsilon
\qquad
\hbox{by Wu [45]},
\cr
& \theta = \textstyle {17\over 41}+\varepsilon
\qquad
\hbox{by Wu [46]},
\cr
& \theta = \textstyle {33\over 80}+\varepsilon
\qquad
\hbox{by Wu [47] and by Zhai [48] (independently)},
\cr
& \theta = \textstyle {40\over 97}+\varepsilon
\qquad
\hbox{by Sargos \& Wu [40]}.
\cr}$$

Inserting our new result on bilinear form ((5.14) of Corollary 9) 
into the argument of [46],
we immediately obtain a slightly better exponent.

\proclaim Proposition 6.
For any $\varepsilon>0$, $x\ge x_0(\ciB, \varepsilon)$ and $y\ge
x^{7/17+\varepsilon}$,
we have
$$\sum_{\scriptstyle x<n\le x+y\atop\scriptstyle b\,\nmid\,n \,
(\forall b\in \cisB)} 1
\gg_{\cisB, \varepsilon} y.$$

Next we shall consider special sets $\ciB$, of the type which occurs in the
applications to modular forms (Theorems 1 and 2).
Let $\P$ be the set of all prime numbers and
$\ciP$ be a subset of $\P$ for which there is a constant $\rho\in [0, 
1]$ such that
$$\big|\ciP\cap [1, x]\big|\ll {x^\rho\over (\log x)^{\Theta_\rho}}
\qquad
(x\ge 2),
\leqno(6.2)$$
where $\Theta_\rho$ is a real constant such that $\Theta_1>1$.
Define
$$\ciB_{\cissP}
:= \ciP \cup \big\{p^2\,\mid\,p\in \P\sset\ciP\big\}
= \{b_k\,\mid\,b_1<b_2<\cdots \,\,\}.$$
Clearly the hypothesis (6.2) guarantees that $\ciB_{\cissP}$
satisfies the condition (6.1).
One can hope to obtain a smaller exponent for this special set of
integers $\ciB_{\cissP}$ than in the general case.
In this direction, Alkan ([1], Theorems 2.2 and 2.3) proved,
by exploiting the structure of the first component $\ciP$ of $\ciB_{\cissP}$,
the following result:
If $y\ge x^{\theta}$ with
$$\theta = \theta(\rho)
= \cases{
{1\over 3} + \varepsilon
& if $\rho = {1\over 2}$,
\cr\noalign{\medskip}
\max\big\{{7\over 19}, \, {23\rho\over 35\rho+16}\big\} + \varepsilon
& if ${1\over 2}<\rho\le 1$,
\cr}
\leqno(6.3)$$
then
$$\sum_{\scriptstyle x<n\le x+y\atop\scriptstyle
b\,\nmid\,n \, (\forall b\in \cisB_{\cisssP})} 1
\gg_{\cissP, \varepsilon} y.
\leqno(6.4)$$
His proof is based on the method of Bantle \& Grupp [3],
using the weight of the form
$$w(n) := \mathop{\sum_{p_1\in {\cal P}_1} \sum_{p_2\in {\cal 
P}_2}}_{p_1p_2\mid n} 1,$$
where
$${\cal P}_i
:= \big\{p\in \P\,\mid\,x^{\delta_i}<p_i\le x^{\delta_i+\varepsilon}\big\}.$$
This leads to estimate a bilinear form of type II:
$$\sum_{m\sim M}\sum_{n\sim N} \varphi_m \psi_n \, r_{mn}(x,y),
\leqno(6.5)$$
where $|\varphi_m|\le 1$, $|\psi_n|\le 1$ and $r_d(x,y)$ is defined in (5.12).
Thus (6.3) is a consequence of the following result of Fouvry \& Iwaniec [15]
with the choice of $x^{\delta_1+\varepsilon} = M$ and
$x^{\delta_2+\varepsilon} = N$:
If $y= x^\theta$,
then for any $\varepsilon>0$ we have
$$\sum_{m\sim M}\sum_{n\sim N} \varphi_m \psi_n \, r_{mn}(x,y)
\ll y x^{-2\varepsilon}
\leqno(6.6)$$
provided
$$\textstyle {7\over 19}<\theta\le {11\over 23},
\qquad
M\le yx^{-\varepsilon'},
\qquad
N\le y^{19/16} x^{-7/16-\varepsilon'}.
\leqno(6.7)$$
It is worth indicating that the condition $M\le yx^{-\varepsilon'}$ 
forces $\delta_1<\theta$,
which obstructs to exploit fully the second component $\big\{p^2\,\mid\,
p\in \P\sset\ciP\big\}$ of $\ciB_{\cissP}$.

In [45] and [46], the third author proposed an improved weighting device,
i.e. replacing ${\cal P}_1$ by a set of quasi-prime numbers ${\cal
M}$ (cf. (7.4) below).
Thanks to the fundamental lemma of sieve ([4], Lemma 4),
we are brought back to estimate the bilinear form of type I defined 
in (5.11). 
Our result (Corollary 9) on bilinear forms of type I has two 
advantages in comparison of (6.7).
Firstly $N$ has a larger range. Secondly there is no condition on $M$ 
as $M\le yx^{-\varepsilon'}$.
The technique of using weights is more effective
if the range of weights can go beyond the natural limit $y$.
In the general case of $\ciB$-free numbers, this is a crucial obstruction.
However the special structure of the second component
$\big\{p^2\,\mid\,p\in \P\sset\ciP\big\}$ of $\ciB_{\cissP}$
allows us to surmount this difficulty
with the result of Filaseta \& Trifonov ([13], (4)).
These two observations and our new estimate for exponential sums 
enable us to improve considerably (6.3) of Alkan.

\proclaim Proposition 7.
Let $0<\rho\le 1$ and $(\kappa, \lambda)$ be an exponent pair.
For any $\varepsilon>0$, $x\ge x_0(\ciP, \varepsilon)$ and $y\ge 
x^{\theta(\rho)}$ with
$$\theta(\rho)
= \max\bigg\{{1\over 3}, \,
{7\rho\over 9\rho + 8}\bigg\} + \varepsilon,
\leqno(6.8)$$
or
$$\theta(\rho)
= \max\bigg\{{\kappa+\lambda\over 1+2\kappa+2\lambda}, \,
{(1+\kappa+2\lambda)\rho\over
(1+2\kappa+2\lambda)\rho+2+2\lambda}\bigg\} + \varepsilon,
\leqno(6.9)$$
we have
$$\sum_{\scriptstyle x<n\le x+y\atop\scriptstyle
b\,\nmid\,n \, (\forall b\in \cisB_{\cisssP})} 1
\gg_{\cissP, \varepsilon} y.$$

When $\rho\le 3(\kappa+\lambda)/(3+2\kappa+2\lambda)$ where $(\kappa, 
\lambda)$ is an exponent pair,
we can obtain a better exponent than that in
Proposition 7.

\proclaim Proposition 8.
For any $\varepsilon>0$, $x\ge x_0(\ciP, \varepsilon)$ and $y\ge 
x^{\theta(\rho)}$ with
$$\theta(\rho)
= \cases{
\displaystyle \max\bigg\{{1\over 4}, \,
{10\rho\over 19\rho+7}\bigg\} + \varepsilon
& if $0\leq \rho<{9\over 17}$,
\cr\noalign{\smallskip}
\displaystyle {3\rho\over 4\rho+3} + \varepsilon
& if ${9\over 17}\le \rho\le 1$,
\cr}
\leqno(6.10)$$
we have
$$\sum_{\scriptstyle x<n\le x+y\atop\scriptstyle
b\,\nmid\,n \, (\forall b\in \cisB_{\cisssP})} 1
\gg_{\cissP, \varepsilon} y.$$

By combining Propositions 7 and 8, we immediately obtain the following result.

\proclaim Corollary 10.
For any $\varepsilon>0$, $x\ge x_0(\ciP, \varepsilon)$ and $y\ge
x^{\theta(\rho)+\eps}$ with 
$$\theta(\rho)
= \cases{
{1\over 4}               & if $0\leq \rho\le {1\over 3}$,
\cr\noalign{\medskip}
{10\rho\over 19\rho + 7}    & if ${1\over 3}<\rho\le
{9\over 17}$,
\cr\noalign{\medskip}
{3\rho\over 4\rho+3}         & if ${9\over 17}<\rho\le
{15\over 28}$,
\cr\noalign{\medskip}
{5\over 16}                 & if ${15\over 28}<\rho\le
{5\over 8}$,
\cr\noalign{\medskip}
{22\rho\over 24\rho + 29}    & if ${5\over 8}<\rho\le
{9\over 10}$,
\cr\noalign{\medskip}
{7\rho\over 9\rho + 8}       & if ${9\over 10}<\rho\le 1$,
\cr}$$
we have
$$\sum_{\scriptstyle x<n\le x+y\atop\scriptstyle
b\,\nmid\,n \, (\forall b\in \cisB_{\cisssP})} 1
\gg_{\cissP, \varepsilon} y.$$

\noindent{\sl Proof}.
The intervals $(0, {1\over 3}]$, $[{1\over 3}, {9\over 17}]$ and
$[{9\over 17}, {15\over 28}]$
come from Proposition 8.
\par
The intervals $[{15\over 28}, {5\over 8}]$ and $[{5\over 8}, {9\over 10}]$
come from (6.9) of Proposition 7 with $(\kappa, \lambda) = ({4\over
18}, {11\over 18})$.
\par
The interval $[{9\over 10}, 1]$
come from (6.8) of Proposition 7.
\hfill
$\square$

\smallskip

{\bf Remark 3.}
(i)
Propositions 7 and 8 improve Alkan's exponent (6.3).
It is worth remarking that we have no restriction $\rho\ge \dm$ as in
[1]. 

(ii)
The parameter $\rho$ can be considered
as a measure of difficulty in the problem of $\ciB_{\cissP}$-free numbers.
Clearly the case $\rho=1$ is the most difficult and $\rho=0$ is the simplest.
In fact when $\ciP$ is empty (so $\rho=0$) the
$\ciB_\emptyset$-free numbers are the square-free integers.
In this case, Filaseta \& Trifonov [13] proved that $\theta={1\over
5}+\eps$ is admissible.
However our method only gives $\theta={1\over 4} + \varepsilon$.
It seems interesting to generalise the method of Filaseta \& Trifonov
to the case of $\ciB_{\cissP}$-free numbers (at least for small values
of $\rho$).

(iii) The function $\theta(\rho)$ is continuous, increasing, and 
$\theta({9\over 17})={9\over 29}$, $\theta({15\over 28})={5\over 16}$,
$\theta({9\over 10})={9\over 23}$, $\theta(1)={7\over 17}$. 

\smallskip

If we relax the multiplicative constraint by removing the square-free
assumption, we can prove a better result for $\rho\le {1\over 3}$.

\proclaim Proposition 9.
Suppose that $0\leq \rho\le 1$.
For any $\varepsilon>0$, $x\ge x_0(\ciP, \varepsilon)$ and $y\ge 
x^{\rho/(1+\rho)+\varepsilon}$, we have
$$\sum_{\scriptstyle x<n\le x+y\atop\scriptstyle
b\,\nmid\,n \, (\forall b\in \cisP)} 1
\gg_{\cissP, \varepsilon} y.$$

\goodbreak
\vskip 5mm

\noindent{\bf \S\ 7. Proof of Proposition 7}

\medskip

We begin by describing our weight function.
Let $\theta$, $\delta_1$ and $\delta_2$ be some parameters such that
$$\textstyle{1\over 4}+\varepsilon\le \theta<{1\over 2},
\quad
\varepsilon<\delta_2+2\varepsilon<\delta_1+\varepsilon<\theta/\rho,
\quad
\delta_1 + \delta_2<1,
\quad
\delta_1 + \delta_2 + \theta/\rho>1.
\leqno(7.1)$$
Introduce two sets
$$\leqalignno{
{\cal M}
& := \bigl\{m\in \N\,\mid\,x^{\delta_1}<m\le x^{\delta_1+\varepsilon}, \,
p\mid m\Rightarrow p\ge x^\eta \big\},
& (7.2)
\cr\noalign{\smallskip}
{\cal P}
& := \big\{p\in \P\,\mid\,x^{\delta_2}<p\le x^{\delta_2+\varepsilon}\big\},
& (7.3)
\cr}$$
where $\eta = \eta(\ciP, \varepsilon)>0$ is a (small) parameter chosen later.

Our weight function is defined by
$$c(n) := \mathop{\sum_{m\in {\cal M}} \sum_{p\in {\cal P}}}_{mp\mid n} 1.
\leqno(7.4)$$
Put
$$A :=
\sum_{\scriptstyle x<n\le x+y\atop\scriptstyle b\;\nmid \, n \,
(\forall b\in \cisB_{\cisssP})}
c(n).
\leqno(7.5)$$
  From (7.1), (7.2) and (7.3), it is easy to see that
$$c(n)\le 2^{1/\eta}/\varepsilon
\qquad(n\le 2x),
\leqno(7.6)$$
which implies
$$\sum_{\scriptstyle x<n\le x+y\atop\scriptstyle b\;\nmid \, n \,
(\forall b\in \cisB_{\cisssP})} 1
\ge \varepsilon 2^{-1/\eta} A.
\leqno(7.7)$$
In order to prove Proposition 7, it is sufficient to show that
$$A\gg_{\cissP, \varepsilon} y.
\leqno(7.8)$$
For this, we let $\ell := \ell(\ciP, \varepsilon)\in \N$ be a
positive  integer such that
$$\sum_{k=\ell+1}^\infty {1\over b_k}<{B_{\cissP}\varepsilon^3\over
\eta 2^{1/\eta+2}},
\leqno(7.9)$$
where
$$B_{\cissP}
:= \prod_{p\in \cisP}\bigg(1- {1\over p}\bigg)
\prod_{p\in \P\sset\cisP}\bigg(1- {1\over p^2}\bigg)$$
is the natural density of the sequence of
$\ciB_{\cissP}$-free numbers.

Clearly we can write
$$A\ge A_1-A_2-A_3
\leqno(7.10)$$
where
$$\eqalign{
A_1 & :=\sum_{\scriptstyle x<n\le x+y\atop\scriptstyle b_k\,\nmid \,
n \, (\forall\, k\le \ell)}
c(n),
\cr
A_2 & :=\sum_{\scriptstyle b_\ell<b\le y\atop\scriptstyle b\in \cisB_{\cisssP}}
\sum_{\scriptstyle x<n\le x+y\atop\scriptstyle b\mid n}c(n),
\cr
A_3 & :=\sum_{\scriptstyle y<b\le x\atop\scriptstyle b\in \cisB_{\cisssP}}
\sum_{\scriptstyle x<n\le x+y\atop\scriptstyle b\mid n} c(n).
\cr}$$

We shall see that $A_2$ and $A_3$ are negligible and
$A_1$ gives the desired principal term.
The required estimates for $A_2$ and $A_3$ will be offered by the
next two lemmas.

\smallskip

\proclaim Lemma 7.1.
We have
$$A_2\le {B_{\cissP}\varepsilon^2\over 2\eta}y.$$

\noindent{\sl Proof}.
By (7.6), it follows that
$$\eqalign{A_2
& \le {2^{1/\eta}\over \varepsilon}
\sum_{\scriptstyle b_\ell<b\le y\atop\scriptstyle b\in \cisB_{\cisssP}}
\sum_{\scriptstyle x<n\le x+y\atop\scriptstyle b\mid n} 1
\cr
& \le {2^{1/\eta}\over \varepsilon}
\sum_{{\scriptstyle b_\ell<b\le y}
\atop{\scriptstyle b\in \cisB_{\cisssP}}}
{2y\over b},
\cr}$$
which implies the required inequality in view of (7.9).
\hfill
$\square$

\smallskip

\proclaim Lemma 7.2.
There is a constant $C(\ciP, \varepsilon)$ such that
$$A_3\le {C(\ciP, \varepsilon) 2^{1/\eta}\over (\log x)^{1/2}} y.$$

\noindent{\sl Proof}.
According to the definition of $\ciB_{\cissP}$, we can write
$$\leqalignno{A_3
& = \sum_{\scriptstyle y<p\le x^{\theta/\rho}(\log x)^{(\Theta_\rho-1/2)/\rho}
\atop\scriptstyle p\in \cisP}
\sum_{\scriptstyle x<n\le x+y\atop\scriptstyle p\mid n} c(n)
& (7.11)
\cr
& \qquad
+ \sum_{\scriptstyle x^{\theta/\rho}(\log x)^{(\Theta_\rho-1/2)/\rho}<p\le x
\atop\scriptstyle p\in \cisP}
\sum_{\scriptstyle x<n\le x+y\atop\scriptstyle p\mid n} c(n)
\cr
& \qquad
+ \sum_{\scriptstyle y<q^2\le y^2\log x\atop\scriptstyle q\in \P\sset\cisP}
\sum_{\scriptstyle x<n\le x+y\atop\scriptstyle q^2\mid n} c(n)
\cr
& \qquad
+ \sum_{\scriptstyle y^2\log x<q^2\le x\atop\scriptstyle q\in \P\sset\cisP}
\sum_{\scriptstyle x<n\le x+y\atop\scriptstyle q^2\mid n} c(n)
\cr
& =: A_{3,1} + A_{3,2} + A_{3,3} + A_{3,4}.
\cr}$$

For $p>y$, there is at most an integer $n\in (x, x+y]$ such that $p\mid n$.
Thus (7.6) and (6.2) imply that
$$\eqalign{A_{3,1}
& \le {2^{1/\eta}\over \varepsilon}
\sum_{\scriptstyle p\le x^{\theta/\rho}(\log x)^{(\Theta_\rho-1/2)/\rho}
\atop\scriptstyle p\in \cisP} 1
\cr
& \ll {2^{1/\eta}\over \varepsilon}
{\big(x^{\theta/\rho}(\log x)^{(\Theta_\rho-1/2)/\rho}\big)^\rho\over
(\log x)^{\Theta_\rho}}
\cr
& \ll {2^{1/\eta}\over \varepsilon (\log x)^{1/2}} y.
\cr}$$

The definition of $c(n)$ allows us to write
$$A_{3,2}
= \sum_{\scriptstyle x^{\theta/\rho}(\log x)^{(\Theta_\rho-1/2)/\rho}<p\le x
\atop\scriptstyle p\in \cisP}
\sum_{m\in {\cal M}} \sum_{p'\in {\cal P}}
\sum_{\scriptstyle x<n\le x+y\atop\scriptstyle p\mid n, \, mp'\mid n} 1.$$
The hypothesis
$\delta_2+2\varepsilon<\delta_1+\varepsilon<\theta/\rho$ and $p\in
\ciP$
imply $(p, mp')=1$.
Thus $pmp'\mid n$.
Since
$$pmp'>x^{\theta/\rho +\delta_1+\delta_2}(\log
x)^{(\Theta_\rho-1/2)/\rho}\ge 2x,$$
the sum over $n$ must be empty. Therefore $A_{3,2}=0$.

We have
$$A_{3,3}
\le {2^{1/\eta}\over \varepsilon}
\sum_{\scriptstyle q\le y(\log x)^{1/2}\atop\scriptstyle q\in \cisP} 1
\ll {2^{1/\eta}\over \varepsilon (\log x)^{1/2}} y.$$

The term $A_{3,4}$ will be treated by the method of Filaseta \&
Trifonov [13].
Defining
$$S(t_1, t_2)
:= \{d\in (t_1, t_2]\,\mid\, \hbox{there is an integer $k$ such that}
\,\, kd^2\in (x, x+y]\},$$
we can deduce, in view of (7.6), that
$$\eqalign{A_{3,4}
& \le \varepsilon^{-1} 2^{1/\eta}
\sum_{\scriptstyle y^2\log x<q^2\le x\atop\scriptstyle q\in \P\sset\cisP}
\sum_{\scriptstyle x<n\le x+y\atop\scriptstyle q^2\mid n} 1
\cr
& \le \varepsilon^{-1} 2^{1/\eta} \big|S\big(y(\log x)^{1/2},
x^{1/2}\big)\big|.
\cr}$$
We split $\big(y(\log x)^{1/2}, x^{1/2}\big]$ into dyadic intervals
$(x^\phi, 2x^\phi]$ and write
$$\eqalign{A_{3,4}
& \le \varepsilon^{-1} 2^{1/\eta} (\log x)
\max_{\theta\le \phi\le 1/2} \big|S(x^\phi, 2x^\phi)\big|.
\cr}$$
According to ([13], (4)), we have
$$\big|S(x^\phi, 2x^\phi)\big|\ll x^{(1-\phi)/3}$$
for $y(\log x)^{1/2}\le x^\phi\le 2x^{1/2}$,
and thus infer with the hypothesis $\theta>{1\over 4}+\varepsilon$ that
$$A_{3,4}
\ll \varepsilon^{-1} 2^{1/\eta} x^{-\varepsilon'} y.$$

Now inserting the estimates for $A_{3,j}$ into (7.11), we obtain the
required inequality.
\hfill
$\square$

\medskip

Next we shall treat the principal term $A_1$.
It is convenient to introduce some notation.
For each $\sigma = \{k_1, \dots, k_i\}\subset \{1, \dots, \ell\}$, we write
$|\sigma|=i$ and $d_\sigma=b_{k_1}b_{k_2}\cdots b_{k_i}$
with the convention $|\emptyset|=0$ and $d_\emptyset=1$,
where $\emptyset$ denotes the empty set.

\smallskip

\proclaim Lemma 7.3.
For $x\ge x_0(\ciP, \varepsilon)$, we have
$$A_1
\ge {B_{\cissP} \varepsilon^2\over \eta} y + R,$$
where
$$R := \sum_{\sigma \subset \{1, \dots, \ell\}}
(-1)^{|\sigma|} \sum_{m\in {\cal M}}\sum_{p\in {\cal P}}
r_{d_\sigma mp}(x,y).
\leqno(7.12)$$

\noindent{\sl Proof}.
Since $(b_j, b_k)=1\,\,(j\not=k)$, we can write
$$\eqalign{A_1
& = \sum_{\sigma \subset \{1, \dots, \ell\}} (-1)^{|\sigma|}
\sum_{\scriptstyle x<n\le x+y\atop\scriptstyle d_\sigma\mid n} c(n)
\cr
& = \sum_{\sigma \subset \{1, \dots, \ell\}}
(-1)^{|\sigma|}\sum_{m\in {\cal M}}\sum_{p\in {\cal P}}
\sum_{\scriptstyle x<n\le x+y\atop\scriptstyle d_\sigma\mid n, \;  mp\mid n} 1.
\cr}$$
Clearly for any $\sigma\subset \{1, \dots, \ell\}$,
any $m\in {\cal M}$ and any $p\in {\cal P}$ with $x\ge x_0(\ciP, \varepsilon)$,
we have $(d_\sigma,mp)=1$ in view of (7.1)--(7.3).
Hence it follows that
$$\leqalignno{A_1
& =\sum_{\sigma \subset \{1, \dots, \ell\}} (-1)^{|\sigma|}\sum_{m\in {\cal M}}
\sum_{p\in {\cal P}}
\sum_{\scriptstyle x<n\le x+y\atop\scriptstyle d_\sigma mp\mid n} 1
& (7.13)
\cr
& = y \sum_{\sigma \subset \{1, \dots, \ell\}}
{(-1)^{|\sigma|}\over d_\sigma}
\sum_{m\in {\cal M}}{1\over m}
\sum_{p\in {\cal P}}{1\over p}+R,
\cr}$$
where $R$ is defined in (7.12).

It is easy to see that
$$\sum_{\sigma \subset \{1, \dots, \ell\}}
{(-1)^{|\sigma|}\over d_\sigma}
= \prod_{k=1}^\ell \bigg(1-{1\over b_k}\bigg)\ge B_{\cissP}
\leqno(7.14)$$
and
$$\sum_{p\in {\cal P}}{1\over p}
= \log\bigg({\delta_2+\varepsilon\over \delta_2}\bigg)
+ O\bigg({1\over \log x}\bigg)
\ge \varepsilon
\leqno(7.15)$$
for $x\ge x_0(\ciP, \varepsilon)$.

In order to estimate the sum over $m$,
we need  the following result of Friedlander ([16], Lemma~2):
{\sl Let $w(t)$ be Buchstab's function
$$w(t)=1/t \quad (1\le t\le 2),
\qquad
\big(t w(t)\big)' = w(t-1) \quad (t\ge 2).$$
Assume $x>1$ and $z = x^{1/t}$ with $t\ge 1$.
Then we have uniformly for $t\ge 2$}
$$\sum_{n\le x, \, p\mid n\Rightarrow p\ge z} 1
= w(t){x\over \log z} + O\bigg({x\over \log^2z}\bigg).$$
  From this, an integration by part deduces that
$$\eqalign{\sum_{m\in {\cal M}} {1\over m}
& = \int_{x^{\delta_1}}^{x^{\delta_1+\varepsilon}}
{1\over t} \d \bigg(\sum_{n\le t, \, p\mid n\Rightarrow p\ge x^\eta} 1\bigg)
\cr
& = {1\over\eta\log x}
\bigg\{
w\bigg({\delta_1+\varepsilon\over\eta}\bigg) - w\bigg({\delta_1\over\eta}\bigg)
\bigg\}
+ O\bigg({1\over\eta^2\log x}\bigg)
+ \int_{\delta_1/\eta}^{(\delta_1+\varepsilon)/\eta} w(u) \d u.
\cr}$$
In view of the well known relation
$$w(t)\rightarrow e^{-\gamma}
\qquad
(t\to \infty),$$
where $\gamma$ is Euler's constant, we immediately see
$${\varepsilon\over 2\eta}
\le \sum_{m\in {\cal M}}{1\over m}
\le {\varepsilon\over \eta}
\leqno(7.16)$$
for $x\ge x_0(\ciP, \varepsilon)$.

Now the expected inequality follows from (7.13)--(7.16).
\hfill
$\square$

\smallskip

The next lemma gives the desired estimate for the error term $R$
defined in (7.12).

\proclaim Lemma 7.4.
Let $s$ be a real number such that
$$s\ge 3
\qquad \hbox{and} \qquad
s\eta<\textstyle{1\over 2}\varepsilon <\textstyle{1\over 4}.
\leqno(7.17)$$
If
$$\cases{
{1\over 3}<\theta\le {5\over 11},
\cr\noalign{\medskip}
\delta_2\le (9\theta-3)/4-\varepsilon',
\cr\noalign{\medskip}
\delta_1+\delta_2\le 1-\varepsilon',
\cr}
\leqno(7.18)$$
or
$$\cases{
(\kappa+\lambda)/(1+2\kappa+2\lambda)<\theta\le
(\kappa+\lambda)/(2\kappa+\lambda),
\cr\noalign{\medskip}
\delta_2\le
[(1+2\kappa+2\lambda)\theta-\kappa-\lambda]/(1+\lambda)-\varepsilon',
\cr\noalign{\medskip}
\delta_1+\delta_2\le 1-\varepsilon',
\cr}
\leqno(7.19)$$
then we have
$$|R|\le C_1(\varepsilon)2^{\ell(\cisP, \varepsilon)}
\big(\eta^{-1} s^{-s} + x^{-\varepsilon/4}\big)y,$$
where $C_1(\varepsilon)$ is a positive constant depending on
$\varepsilon$ only.

\noindent{\sl Proof.}
For each $\sigma\subset\lbrace 1, \dots,\ell\rbrace$, we define
$$R(\sigma)
:= \sum_{m\in {\cal M}}\sum_{p\in {\cal P}} r_{d_{\sigma}mp}(x,y).$$
We shall transform $R(\sigma)$ into a bilinear form of type I by using
the fundamental lemma of sieve ([4], Lemma 4):
{\sl Let $z=x^\eta$ and $Q=z^s$ with $s\ge 3$.
There are two sequences $\{\lambda^\pm_q\}_{q\le Q}$
such that}
$$\leqalignno{
& \,\,\, |\lambda^\pm_q|\le 1,
\qquad
\lambda^\pm_q=0
\quad(q>Q),
& (7.20)
\cr\noalign{\medskip}
& \cases{
(\lambda^-*{\bf 1})(n) = (\lambda^+*{\bf 1})(n)=1
& if $p\mid n \Rightarrow p\ge z$,
\cr\noalign{\medskip}
(\lambda^-*{\bf 1})(n)\le 0\le (\lambda^+*{\bf 1})(n)
& otherwise,
\cr}
& (7.21)
\cr\noalign{\smallskip}
& \,\, \sum_{q\le Q}{{\lambda}^{\pm}_q\over q}
= \{1 + O(s^{-s})\}
\prod_{p<z}\bigg(1-{1\over p}\bigg).
& (7.22)
\cr}$$

With the help of (7.21), we can write
$$\leqalignno{R(\sigma)
& = \sum_{m\in {\cal M}}\sum_{p\in {\cal P}}
\sum_{\scriptstyle x<n\le x+y\atop\scriptstyle d_{\sigma}mp\,\mid\,n} 1
- \sum_{m\in {\cal M}}\sum_{p\in {\cal P}} {y\over d_{\sigma}mp}
& (7.23)
\cr
& \le \sum_{x^{\delta_1}<m\le x^{\delta_1+\varepsilon}} (\lambda^+*{\bf 1})(m)
\sum_{p\in {\cal P}}
\sum_{\scriptstyle x<n\le x+y\atop\scriptstyle d_{\sigma}mp\,\mid\,n} 1
\cr
& \qquad
- \sum_{x^{\delta_1}<m\le x^{\delta_1+\varepsilon}} (\lambda^-*{\bf 1})(m)
\sum_{p\in {\cal P}} {y\over d_{\sigma}mp}
\cr
& = \sum_{x^{\delta_1}<m\le x^{\delta_1+\varepsilon}} (\lambda^+*{\bf 1})(m)
\sum_{p\in {\cal P}}
r_{d_{\sigma}mp}(x, y)
\cr
& \qquad
+ \sum_{x^{\delta_1}<m\le x^{\delta_1+\varepsilon}}
[(\lambda^+ - \lambda^-)*{\bf 1}](m)
\sum_{p\in {\cal P}}
{y\over d_{\sigma}mp}
\cr
& =:  R_1(\sigma)+R_2(\sigma).
\cr}$$

Clearly
$$\eqalign{|R_2(\sigma)|
& \le {y\over d_{\sigma}}
\sum_{q\le Q} \bigg|{\lambda_q^+ - \lambda_q^-\over q}\bigg|
\sum_{x^{\delta_1}/q<m\le x^{\delta_1+\varepsilon}/q} {1\over m}
\sum_{p\in {\cal P}} {1\over p}
\cr
& \le {y\over d_\sigma}
\sum_{q\le Q} \bigg|{\lambda_q^+ - \lambda_q^-\over q}\bigg|
\cdot \big\{\varepsilon \log x + O(qx^{-\delta_1})\big\} \cdot 2\varepsilon
\cr}$$
On the other hand, (7.22) implies that
$$\sum_{q\le Q} \bigg|{\lambda_q^+ - \lambda_q^-\over q}\bigg|
\ll {s^{-s}\over \log z}
\ll {s^{-s}\over \eta \log x}.$$
Inserting it into the preceding estimate, we find that
$$\big|R_2(\sigma)\big|
\ll \big(\eta^{-1}s^{-s} + Qx^{-\delta_1}\big)y.
\leqno(7.24)$$

It remains to estimate $R_1(\sigma)$.
Let $\psi_n$ be the characteristic function of the set ${\cal P}$.
Since
$$r_{d_\sigma qmn}(x,y) = r_{mn}\bigg({x\over d_\sigma q}, {y\over
d_\sigma q}\bigg),$$
we can write
$$R_1(\sigma)
=\sum_{q\le Q}\lambda_q^+
\sum_{x^{\delta_1}/q<m\le x^{\delta_1+\varepsilon}/q} \,
\sum_{x^{\delta_2}<n\le x^{\delta_2+\varepsilon}} \psi_n
r_{mn}\bigg({x\over d_\sigma q}, {y\over d_\sigma q}\bigg).$$
We split $(x^{\delta_1}/q, \, x^{\delta_1+\varepsilon}/q]$
and $(x^{\delta_2}, \, x^{\delta_2+\varepsilon}]$ into dyadic intervals
$(M, 2M]$ and $(N, 2N]$, respectively.
In view of (7.22), we have for $x\ge x_0(\ciP, \varepsilon)$
$$1\le q\le Q=x^{s\eta}<x^{\varepsilon/2}
\qquad\hbox{and}\qquad
1\le d_\sigma<x^{\varepsilon/2}.$$
The hypothesis (7.18) and (7.19) imply that
$$N<\bigg({x\over d_\sigma q}\bigg)^{(9\theta-3)/4-\varepsilon'}
\qquad{\rm and}\qquad
MN\le\bigg({x\over d_\sigma q}\bigg)^{1-\varepsilon'}$$
and
$$N<\bigg({x\over d_\sigma
q}\bigg)^{[(1+2\kappa+2\lambda)\theta-\kappa-\lambda]/(1+\lambda)
-\varepsilon'}
\qquad{\rm and}\qquad
MN\le\bigg({x\over d_\sigma q}\bigg)^{1-\varepsilon'},$$
respectively.
Thus Corollary 9 allows us to deduce that
$$\eqalign{\sum_{x^{\delta_1}/q<m\le x^{\delta_1+\varepsilon}/q}
\sum_{x^{\delta_2}<n\le x^{\delta_2+\varepsilon}} \psi_n
r_{mn}\bigg({x\over d_\sigma q}, {y\over d_\sigma q}\bigg)
& \ll_\varepsilon {y\over d_\sigma q}
\bigg({x\over d_\sigma q}\bigg)^{-\varepsilon}(\log x)^2
\cr
& \ll_\varepsilon x^{-\varepsilon/2} y.
\cr}$$
This estimate and (7.20) imply that
$$R_1(\sigma)\ll_\varepsilon Qx^{-\varepsilon/2}y
\ll_\varepsilon x^{-\varepsilon/4}y.
\leqno(7.25)$$
Combining (7.24) and (7.25),
there is a positive constant $C_1(\varepsilon)>0$ depending on
$\varepsilon$ such that
$$R(\sigma)\le C_1(\varepsilon)
\big(\eta^{-1}s^{-s}+x^{-\varepsilon /4}\big)y.$$

\goodbreak

Similarly we can prove that
$$R(\sigma)\ge -C_1(\varepsilon)
\big(\eta^{-1}s^{-s}+x^{-\varepsilon /2}\bigr)y.$$
Thus
$$|R|\le \sum_{\sigma\subset \{1, \dots, \ell\}}
|R(\sigma)|\le C_1(\varepsilon)2^{\ell(\cisP, \varepsilon)}
\big(\eta^{-1}s^{-s}+x^{-\varepsilon /2}\bigr)y.$$
This completes the proof of Lemma 7.4.
\hfill
$\square$

\medskip

Now we are ready to complete the proof of Proposition 7.

Without loss of generality, we can assume that
$$C_1(\varepsilon) 2^{\ell(\cisP, \varepsilon)}
>   8/B_{\cissP}\varepsilon^2
>16.$$
Take
$$\displaylines{
\eta^{-1} =
\min\big\{{\textstyle {1\over 5}}\varepsilon^{-2}, \,
C_1(\varepsilon) 2^{\ell(\cisP, \varepsilon)}\big\},
\qquad
s=\eta^{-1/2},
\cr\noalign{\smallskip}
\theta
= \textstyle \max\big\{{1\over 3}, \,
{7\rho\over 9\rho + 8}\big\} + \varepsilon'
\quad{\rm or}\quad
\theta
= \max\big\{{\kappa+\lambda\over 1+2\kappa+2\lambda}, \,
{(1+\kappa+2\lambda)\rho\over
(1+2\kappa+2\lambda)\rho+2+2\lambda}\big\} + \varepsilon',
\cr\noalign{\smallskip}
\delta_1 = \theta/\rho - \varepsilon',
\qquad
\delta_2 = 1 - 2\theta/\rho + \varepsilon'.
\cr}$$
It is easy to verify that these choices satisfy the conditions (7.1),
(7.17), (7.18) or (7.19).
Thus Lemmas 7.1--7.5 imply that
$$\eqalign{A
& \ge \bigg(
{B_{\cissP} \varepsilon^2\over 2\eta}
- C_1(\varepsilon)2^{\ell(\cisP, \varepsilon)}
\big(\eta^{-1} s^{-s} + x^{-\varepsilon/4}\big)
- {C(\ciP, \varepsilon)2^{1/\eta}\over (\log x)^{1/2}}
\bigg)
y
\cr
& \gg_{\cissP, \varepsilon} y
\cr}$$
for $x\ge x_0(\ciP, \varepsilon)$.
This completes the proof of (7.8) and hence Proposition 7.
\hfill
$\square$

\vskip 5mm

\noindent{\bf \S\ 8. Proof of Proposition 8}

\medskip

The proof is very similar to that of Proposition 7
so we shall mention only the important points.
As before let $\theta$ and $\delta$ be two parameters such that
$$\textstyle{1\over 4}+\varepsilon\le \theta<{1\over 2},
\qquad
\theta<\delta+2\varepsilon<\min\{\theta/\rho, \, 1\},
\qquad
\delta + \theta/\rho>1.
\leqno(8.1)$$
Let $\eta=\eta(\ciP, \varepsilon)>0$ be a (small) parameter determined later.
Introduce the set
$${\cal M}'
:= \bigl\{m\in \N\,\mid\,x^{\delta}<m\le x^{\delta+\varepsilon}, \,
p\mid m\Rightarrow p\ge x^\eta \big\}.$$
Our weight function is defined to be
$$c'(n) := \sum_{\scriptstyle m\in {\cal M}'\atop\scriptstyle m\mid n} 1$$
and the corresponding weighted sum is
$$A' :=
\sum_{\scriptstyle x<n\le x+y\atop\scriptstyle b\;\nmid \, n \,
(\forall b\in \cisB_{\cisssP})}
c'(n).$$

It is easy to see that
$$c'(n)\le 2^{1/\eta}
\qquad(n\le 2x)
\leqno(8.2)$$
and
$$\sum_{\scriptstyle x<n\le x+y\atop\scriptstyle b\;\nmid \, n \,
(\forall b\in \cisB_{\cisssP})} 1
\ge 2^{-1/\eta} A'.
\leqno(8.3)$$

Let $\ell := \ell(\ciP, \varepsilon)\in \N$ be a positive  integer such that
$$\sum_{k=\ell+1}^\infty {1\over b_k}<{B_{\cissP}\varepsilon^2\over
\eta 2^{1/\eta+2}}.
\leqno(8.4)$$
We can write
$$A'\ge A'_1-A'_2-A'_3
\leqno(8.5)$$
where
$$\eqalign{
A'_1 & :=\sum_{\scriptstyle x<n\le x+y\atop\scriptstyle b_k\,\nmid \,
n \, (\forall\, k\le \ell)}
c'(n),
\cr
A'_2 & :=\sum_{\scriptstyle b_\ell<b\le y\atop\scriptstyle b\in
\cisB_{\cisssP}}
\sum_{\scriptstyle x<n\le x+y\atop\scriptstyle b\mid n} c'(n),
\cr
A'_3 & :=\sum_{\scriptstyle y<b\le x\atop\scriptstyle b\in \cisB_{\cisssP}}
\sum_{\scriptstyle x<n\le x+y\atop\scriptstyle b\mid n} c'(n).
\cr}$$

Similar to Lemmas 7.1, 7.2 and 7.3, we have, for $x\ge x_0(\ciP, \varepsilon)$,
$$\leqalignno{
& A'_2\le {B_{\cissP}\varepsilon^2\over 2\eta}y,
& (8.6)
\cr
& A'_3\le {C 2^{1/\eta}\over (\log x)^{1/2}} y,
& (8.7)
\cr
& A'_1
\ge {B_{\cissP} \varepsilon^2\over \eta} y + R,
& (8.8)
\cr}$$
where
$$R' := \sum_{\sigma \subset \{1, \dots, \ell\}}
(-1)^{|\sigma|} \sum_{m\in {\cal M}'} r_{d_\sigma m}(x,y).$$

Similar to Lemma 7.4, we can prove, by using (5.17) and (5.18) of 
Lemma 5.1 instead of Corollary 9,
that there is a positive constant $C'_1(\varepsilon)$ depending on
$\varepsilon$ only such that
$$|R'|\le C'_1(\varepsilon)2^{\ell(\cisP, \varepsilon)}
\big(\eta^{-1} s^{-s} + x^{-\varepsilon/4}\big)y
\leqno(8.9)$$
provided
$$s\ge 3,
\qquad
s\eta<\textstyle{1\over 2}\varepsilon <\textstyle{1\over 4}
\leqno(8.10)$$
and
$$\cases{
{1\over 4}<\theta<{9\over 29},
\cr\noalign{\medskip}
\delta\le (19\theta-3)/7-\varepsilon',
\cr}
\qquad{\rm or}\qquad
\cases{
{9\over 29}<\theta<{1\over 2},
\cr\noalign{\medskip}
\delta\le 4\theta/3-\varepsilon'.
\cr}
\leqno(8.11)$$

Now take
$$\eta^{-1} =
\min\big\{{\textstyle {1\over 5}}\varepsilon^{-2}, \,
C_1(\varepsilon) 2^{\ell(\cisP, \varepsilon)}\big\},
\qquad
s=\eta^{-1/2}$$
and
$$\cases{\theta
= \max\big\{{1\over 4}, \, {10\rho\over 19\rho+7}\big\} + \varepsilon',
\cr\noalign{\medskip}
\delta = {19\theta-3\over 7} - \varepsilon',
\cr}
\qquad{\rm or}\qquad
\cases{\theta
= {3\rho\over 4\rho+3} + \varepsilon',
\cr\noalign{\medskip}
\delta = {4\theta\over 3} - \varepsilon'.
\cr}$$
It is straightforward to verify that these choices satisfy the 
conditions (8.1),
(8.10) and (8.11).
Thus the relations (8.5)--(8.9) imply
$$A'
\gg_{\cissP, \varepsilon} y$$
for $x\ge x_0(\ciP, \varepsilon)$.
This completes the proof of Proposition 4.
\hfill
$\square$

\vskip 5mm

\noindent{\bf \S\ 9. Proof of Proposition 9}

\medskip

The proof of Proposition 9 (which can in fact be properly described as
a sieve problem in the usual sense) is much simpler than that of
Proposition 8. So we shall mention only the important points.
Let $\theta = \rho/(1+\rho)+2\varepsilon$ and
$${\cal P}''
:= \bigl\{p\in \P\,\mid\,x^{\theta-2\varepsilon}<p\le 
x^{\theta-\varepsilon}\big\}.$$
Define the weight function
$$c''(n) := \sum_{\scriptstyle p\in {\cal P}''\atop\scriptstyle p\mid n} 1$$
and consider the corresponding weighted sum
$$A'' :=
\sum_{\scriptstyle x<n\le x+y\atop\scriptstyle b\;\nmid \, n \, 
(\forall b\in \cisP)}
c'(n).$$
Similarly we can write
$$A''\ge A''_1-A''_2-A''_3,$$
where $A''_j$ is defined as $A'_j$ (replacing $\ciB_{\cisP}$ by $\ciP$).
Now $A''_3$ is easier to treat
(without the corresponding parts $A_{3,3}$ and $A_{3,4}$, see (7.11)).
In view of $\theta-2\varepsilon + \theta/\rho>1$,
we can prove the same estimates for $A''_2$ and $A''_3$.
The error term $R''$, which comes from $A''_1$, can be controlled 
trivially as follows:
$$|R''|
\le \sum_{\sigma \subset \{1, \dots, \ell\}} \sum_{p\in {\cal P}''} |r_{d_\sigma 
p}(x,y)|
\ll_{\cissP} yx^{-\varepsilon}.$$
This completes the proof of Proposition 9.
\hfill
$\square$

\centerline{\bf References}

\medskip

\item{[1]}
{\author E. Alkan},
Nonvanishing of Fourier coefficients of modular forms,
{\it Proc. Amer. Math. Soc.} {\bf 131} (2002), 1673--1680.

\item{[2]}
{\author A. Balog \& K. Ono},
The Chebotarev density theorem in short intervals and some questions of Serre,
{\it J. Number theory} {\bf 91} (2001), 356--371.

\item{[3]}
{\author G. Bantle \& F. Grupp},
On a problem of Erd\H os and Szemer\'edi,
{\it J. Number Theory} {\bf 22} (1986), 280--288.

\item{[4]}
{\author E. Bombieri, J.B. Friedlander \& H. Iwaniec},
Primes in arithmetic progressions to large moduli,
{\it Acta Math}. {\bf 156} (1986), 203--251.

\item{[5]}
{\author A. Booker},
 A test for identifying Fourier coefficients of automorphic forms and
application to Kloosterman sums,
{\it Experiment. Math.}  {\bf 9} (2000), 571--581.

\item{[6]}
{\author J.W.S Cassels},
Local fields, Cambridge Univ. Press 1986.

\item{[7]}
{\author W. Duke \& H. Iwaniec},
Bilinear forms in the Fourier coefficients of half-integral weight
cusp forms and sums over primes, 
{\it Math. Ann.} {\bf 286} (1990), 783--802.

\item {[8]}
{\author W. Duke, J. Friedlander \& H. Iwaniec},
Bounds for automorphic $L$-functions, II,
{\it Invent. math.} {\bf 115} (1994), 219--239.

\item{[9]}
{\author W. Duke \& E. Kowalski},
A problem of Linnik for elliptic curves and mean-value estimates for
automorphic representations,
{\it Invent. math.} {\bf 139} (2000), 1--39.

\item{[10]}
{\author N. Elkies},
Distribution of supersingular primes,
{\it Ast\'erisque} (1992), 127--132.

\item{[11]}
{\author P. Erd\H os},
On the difference of consecutive terms of sequences,
defined by divisibility properties,
{\it Acta Arith}. {\bf 12} (1966), 175--182.

\item{[12]}
{\author J-H Evertse, H. Schlickewei \& W. Schmidt},
Linear equations in variables which lie in a multiplicative group,
{\it Annals of Math.} {\bf 155} (2002), 807--836.

\item{[13]}
{\author M. Filaseta \& O. Trifonov},
On gaps between squarefree numbers II,
{\it J. London Math. Soc.} (2) {\bf 45} (1992), 215--221.

\item{[14]}
{\author B. Fisher}
Distinctness of Kloosterman sums,
in {\it $p$-adic methods in number theory and algebraic geometry},
Contemporary Math. 133, A.M.S 1992, 81--102.

\item{[15]}
{\author E. Fouvry \& H. Iwaniec},
Exponential sums with monomials,
{\it J. Number Theory} {\bf 33} (1989), 311--333.

\item{[16]}
{\author J.B. Friedlander},
Integers free from large and small primes,
{\it Proc. London Math. Soc}. (3) {\bf 33} (1976), 565--576.

\item{[17]}
{\author W. Fulton \& J. Harris},
Representation Theory, a first course,
Grad. Texts in Math. 129, Springer Verlag, 1991.

\item{[18]}
{\author S. Gelbart \& H. Jacquet},
A relation between automorphic representations of $GL(2)$ and
$GL(3)$,
{\it Ann. Sci. \'Ecole Norm. Sup. (4)} {\bf 11} (1978), 471--542.

\item{[19]}
{\author A. Good},
Cusp forms and eigenfunctions of the laplacian,
{\it Math. Ann.} {\bf 255} (1981), 523--548.

\item{[20]}
{\author S.W. Graham \& G. Kolesnik},
{\it Van der Corput's Method of Exponential sums},
Cambridge University Press, Cambridge, 1991.

\item{[21]}
{\author D.R. Heath-Brown},
{\it Mean values of the zeta-function and divisor problems},
in {\it Recent progress in analytic number theory}, Vol. 1 (Durham,
1979), pp. 115--119,
Academic Press, London-New York, 1981.

\item{[22]}
{\author J. Hoffstein \& P. Lockhart},
Coefficients of Maass forms and the Siegel zero,
{\it Annals of Math.} {\bf 140} (1994), 161--181.

\item{[23]}
{\author H. Iwaniec},
{\it Topics in Classical Automorphic Forms},
Graduate Studies in Mathematics, vol. 17, American Mathematical 
Society, Providence, Rhode Island, 1997.

\item{[24]}
{\author H. Iwaniec \& E. Kowalski},
Analytic Number Theory, Colloquium Publications 53, American
Mathematical Society, Providence, Rhode Island, 2004.

\item{[25]}
{\author N. Katz},
Sommes exponentielles,
{\it Ast\'erisque} {\bf 79}, SMF 1980.

\item{[26]}
{\author H. Kim \& P. Sarnak},
Refined estimates towards the Ramanujan and Selberg conjectures,
{\it J. American Math. Soc.} {\bf 16} (2003), 175--181.

\item{[27]}
{\author M.I. Knopp \& J. Lehner},
Gaps in the Fourier series of automorphic forms, 
Analytic Number Theory (Philadelphia 1980),
{\it Lecture Notes in Math.} {\bf 899}, Springer Verlag 1981, 360--381.

\item{[28]}
{\author S. Lang \& H. Trotter},
Frobenius distribution in $GL(2)$ extensions,
Lecture Notes in Mathematics {\bf 504}, Springer-Verlag, Berlin-New 
York, 1976. iii+274 pp.

\item{[29]}
{\author D.H. Lehmer},
Some functions of Ramanujan,
{\it Math. Student} {\bf 27} (1959), 105--116.

\item{[30]}
{\author H.-Q. Liu \& J. Wu},
Numbers with a large prime factor,
{\it Acta Arith}. {\bf 89} (1999), 163--187.

\item{[31]}
{\author C. Moeglin \& J-L. Waldspurger},
P\^oles des fonctions $L$ de paires pour $GL(N)$, appendix to
Le spectre r\'esiduel de $GL(N)$,
{\it Ann. Sci. \'Ecole Norm. Sup. (4)} {\bf 22} (1989), 605--674.

\item{[32]}
{\author G. Molteni},
Upper and lower bounds at $s=1$ for certain Dirichlet series
with Euler product,
{\it Duke Math. J.} {\bf 111} (2002), 133--158.

\item{[33]}
{\author V.K. Murty},
Frobenius distributions and Galois representations,
in: {\it Automorphic forms, automorphic representations and arithmetic}
(Fort Worth, TX, 1996), 193--211, Proc. Sympos. Pure Math. {\bf 66}, Part 1,
Amer. Math. Soc., Providence, RI, 1999.

\item{[34]}
{\author Y. Petridis},
On squares of eigenfunctions for the hyperbolic plane and a new bound
on certain $L$-series, 
{\it IMRN} {\bf 1995}, n. 3, 111--127. 

\item{[35]}
{\author V.A. Plaksin},
Distribution de ${\cal B}$-free numbers,
{\it Matematicheskie Zametki} {\bf 47} (1990), 69--77.

\item{[36]}
{\author R.A Rankin},
Contributions to the theory of Ramanujan's function $\tau(n)$ and 
similar arithmetical functions II.
The order of the Fourier coefficients of integer modular forms,
{\it Proc. Cambridge Phil. Soc.} {\bf 35} (1939), 357--372.

\item{[37]}
{\author K.A. Ribet},
Galois representations attached to eigenforms with Nebentypus,
in {\it Modular functions of one variable V}
(Proc. Second Internat. Conf., Univ. Bonn, Bonn, 1976),
Lecture Notes in Math. {\bf 601} (1977), 17--51, Springer \& Berlin.

\item{[38]}
{\author O. Robert \& P. Sargos},
A third derivative test for mean values of exponential sums
with application to lattice point problems,
{\it Acta Arith.} {\bf 106} (2003), 27--39.

\item{[39]}
{\author O. Robert \& P. Sargos},
Three dimensional  exponential sums with monomials,
{\it J. Reine angew Math.}, to appear.

\item{[40]}
{\author P. Sargos \& J. Wu},
Multiple exponential sums with monomials and their applications in
number theory,
{\it Acta Math. Hungar.} {\bf 87 (4)} (2000), 333--354.

\item{[41]}
{\author P. Sarnak},
Integral of products of eigenfunctions,
{\it IMRN} {\bf 1994}, n. 6, 251--260.

\item{[42]}
{\author A. Selberg},
Bemerkungen \"uber eine Dirichletsche Reihe, die mit der Theorie der 
Modulformen nahe verbunden ist,
{\it Arch. Math. Naturvid.} {\bf 43} (1940), 47--50.

\item{[43]}
{\author J.-P. Serre},
Quelques applications du th\'eor\`eme de densit\'e de Chebotarev,
{\it Inst. Hautes \'Etudes Sci. Publ. Math.} {\bf 54} (1981), 323--401.

\item{[44]}
{\author E. Szemer\'edi},
On the difference of consecutive
terms of sequences, defined by divisibility properties II,
{\it Acta Arith}. {\bf 23} (1973), 359--361.

\item{[45]}
{\author J. Wu},
{\it Sur trois questions  classiques de crible} :
{\it nombres premiers jumeaux, nombres $P_2$ et nombres ${\cal B}$--libres},
Th\`ese de Doctorat, Universit\'e Paris-Sud, 1990.

\item{[46]}
{\author J. Wu},
Nombres ${\cal B}$--libres dans les petits intervalles,
{\it Acta Arith.} {\bf 65} (1993), 97--116.

\item{[47]}
{\author J. Wu},
Distance entre nombres ${\cal B}$-libres cons\'ecutifs,
{\it Pr\'epublications de l'Institut \'Elie Cartan} N$^\circ$ {\bf 32} (1994),
Universit\'e Henri Poincar\'e (Nancy 1).

\item{[48]}
{\author W.G. Zhai} ,
Number of ${\cal B}$-free numbers in short intervals,
{\it Chinese Sci. Bull.} {\bf 45} (2000), no. 3, 208--212.


\vskip 10mm

{\author E. Kowalski} :
{Universit\'e Bordeaux I - A2X},
{351, cours de la Lib\'eration},
{33405 Talence Cedex},
{France}
\par
{e-mail : emmanuel.kowalski@math.u-bordeaux1.fr}
\bigskip
{\author O. Robert \& J. Wu} :
{Institut Elie Cartan},
{UMR 7502 UHP-CNRS-INRIA},
{Universit\'e Henri Poincar\'e (Nancy 1)},
{54506 Vand\oe uvre-l\`es-Nancy},
{France}
\par
{e-mail: robert@iecn.u-nancy.fr\quad wujie@iecn.u-nancy.fr}
\end